\documentclass[12pt]{article}
\usepackage{amsmath}
\usepackage{graphicx,psfrag,epsf}
\usepackage{enumerate}
\usepackage{natbib}
\usepackage{amssymb}
\usepackage[colorlinks,citecolor=blue,urlcolor=blue]{hyperref}
\usepackage{mathrsfs}
\usepackage{amsfonts,anysize,amsmath,indentfirst,geometry,xcolor}
\usepackage{bm}
\usepackage{amsthm}
\usepackage{multirow}
\usepackage{chngpage}
\usepackage{array}
\usepackage{rotating}
\usepackage{setspace}
\usepackage{graphicx}
\usepackage{algorithm}
\usepackage{algorithmic}
\usepackage{epsfig}
\usepackage{amsbsy}
\usepackage{bbm}

\newtheorem{thm}{\underline{\bf Theorem}}

\newtheorem{proposition}{Proposition}

\newtheorem{remark}{\underline{\bf Remark}}

\def\boxit#1{\vbox{\hrule\hbox{\vrule\kern6pt\vbox{\kern6pt#1\kern6pt}\kern6pt\vrule}\hrule}}
\def\bqe{\begin{eqnarray}}
\def\eqe{\end{eqnarray}}

\newcommand{\blind}{0}

\begin{document}
\doublespacing

\def\spacingset#1{\renewcommand{\baselinestretch}%
{#1}\small\normalsize} \spacingset{1}


\if0\blind
{
  \title{\bf Rank-based Model Average Estimator for Average Treatment Effects}
  \author{ Jing Lv\\
    School of Mathematics and Statistics,\\
     Southwest University, China\\
     Jialiang Li \\
         Department of Statistics \& Data Science,\\
        National University of Singapore, Singapore\\
    }
  \maketitle
} \fi

\if1\blind
{
  \bigskip
  \bigskip
  \bigskip
  \begin{center}
    {\LARGE\bf Rank-based Model Average Estimator for Average Treatment Effects}
\end{center}
  \medskip
} \fi

\bigskip
\begin{abstract}
\spacingset{1.85}
Model averaging emerges as an attractive ensemble technique that strategically synthesizes all available information from potential candidate models to enhance estimation accuracy. Despite of having been widely practiced in various areas in biostatistics and economics, its application to causal inference is rather limited. We propose to estimate the average treatment effect (ATE) by forming the weighted combination of the
robust rank-based estimators obtained from a series of candidate models, where the model weights, selected via the growing window cross-validation, dynamically adjust the importance of candidate models. The proposed robust rank-based model averaging ({\sf RRMA}) is not only robust to handle potential outliers and heavy-tailed error distributions but also capable to reduce the risk of model mis-specification and enhance the flexibility and the accuracy of estimation. Theoretically, we establish the desirable asymptotic properties of {\sf RRMA}, including the asymptotic optimality, consistency and asymptotic distribution. Moreover, motivated by the asymptotic distribution of {\sf RRMA}, an easy-to-implement confidence interval is also developed for the ATE. Simulation experiments and an empirical application demonstrate the superiority of {\sf RRMA} over existing competitors in estimation accuracy, coverage probabilities and widths of confidence intervals.
\end{abstract}

\\

\spacingset{1.85}
\section{Introduction}\label{sec:intro}

We consider the problem of evaluating the average treatment effects (ATE) caused by a policy intervention on a single unit in panel data settings with a large number of control units and with many pretreatment periods. To assess such a treatment effect, one often needs to measure the difference between the outcomes of the treated unit under the treatment of interest and the outcomes of the same unit without receiving the treatment. However, we cannot observe both potential outcomes simultaneously.

The synthetic control method (SCM) was first introduced by \cite{AG03} and quickly became a popular approach for estimating the counterfactual outcomes. The basic idea of the SCM is to construct a weighted average of estimated outcomes from control units,
where the weights are often selected by minimizing the squared imbalance between the pretreatment outcomes for the treated unit and the weighted control units. In particular, \cite{ADH10} applied the SCM to investigate whether California's Proposition 99 (a large-scale tobacco control program implemented in California in 1988) significantly impacts cigarette sales in California.
\cite{L20} derived the asymptotic distributions of the modified synthetic control estimator under both trend-stationary and unit-root cases, and also constructed confidence intervals by using a subsampling method.
\cite{B21} studied the asymptotic properties of the synthetic control estimators when both the number of pretreatment
periods and control units tend to infinity. When there is a large discrepancy between the values of the predictors for the treated unit and the values of the predictors for the synthetic control unit, the original SCM may be biased. To overcome this issue, \cite{BFR21} introduced the Augmented SCM that utilizes ridge regression to adjust for mismatches in the pretreatment outcomes. Other related literature for the SCM include \cite{LB17,CWZ21} and among others.

In practice, when a large number of control units are available for predicting the counterfactual outcomes, one usually needs to decide which control unit is optimal. To address this issue, \cite{HCW12} suggested to select the best model according to the Akaike
information criterion (AIC) or the corrected AIC (CAIC). \cite{LB17} recommended the LASSO method to select important control units, resulting in more accurate out-of-sample
prediction results.
Further, \cite{MM21} developed a modified LASSO estimator for the counterfactual outcomes under the framework of high-dimensional, nonstationary data with either deterministic and/or stochastic trends.
However, using a single selected model may overlook some useful information contained in other models and ignores the uncertainty across different candidate models. Instead of selecting a specific ``best" model based on a
criterion (i.e., AIC or CAIC) or regularization strategies (i.e., LASSO), model averaging aims to combine estimators through constructing a weighted ensemble and leads to a lower risk of model mis-specification than single-model methods. Over the past decades, there is a large and growing body of literature on frequentist model
averaging, see \cite{hr12,li2015flexible,HMZZ22,FYT23,ZL23,JLLC24}. However, so far, model averaging for causal inference is rather limited. Recently, \cite{GLW15} proposed to use a model average to improve the estimation of average treatment effects but does not provide any theoretical justification. Combined with the partition and matching technique, \cite{SZZ24} introduced a jackknife model averaging estimator for the conditional average treatment effects under the heteroskedastic error settings.

All the aforementioned works \citep{GLW15,LB17,L20,SZZ24} on estimation of the ATE rely on Gaussian or sub-Gaussian error assumptions.
When the error distribution deviates from normality, these estimators not only exhibit pronounced sensitivity to outliers but also suffer a huge loss in efficiency.  Wilcoxon rank regression (or rank regression for short) is a distribution-free, robust and efficient modeling tool.
\cite{WL09} proposed a robust shrinkage-type variable selection procedure by a combination of rank loss function and SCAD penalty in the low dimensional setting. Further, \cite{WPBLW20} proposed a new ${\rm L}_1$
regularized estimator based on rank loss function with theoretical guarantees under the ultrahigh-dimensional setting. More related literature for rank-based methodology can refer to \cite{BW05,ZWZ24}. From our review, few studies applied rank regression to statistical inference for the ATE.

Our primary goal is to introduce a robust model averaging procedure
for estimating the ATE in this paper. Specifically, our contribution in this paper is threefold.
First, we propose a combined estimator for the ATE by taking the weighted average of the robust rank-based estimators obtained from a series of candidate models, and the model weights selected by the cross validation automatically. Compared to the existing SCM based on a single model, our method markedly reduces the risk of model mis-specification and hence generally provides more flexible and accurate predictive inference for the ATE.
Second, we systematically study the theoretical properties.
Specifically, when all candidate models are mis-specified, the asymptotic optimality and consistency of the model weight estimator, as well as the asymptotic distribution of {\sf RRMA}, are established. In the other scenario where there is at least one quasi-correct model in the candidate model set, we prove that the summation of weights assigned to the quasi-correct candidate submodels tends to one in probability, and also derive {\sf RRMA}'s consistency and asymptotic distribution.
Third, we design a new algorithm to build a series of non-nested candidate
models by a combination of the growing window cross-validation and the average absolute
prediction error. This algorithm works very well to balance prediction accuracy and
model complexity.
To the best of our knowledge, there are no data-driven approaches to design the optimal non-nested candidate models in the field of causal inference. Our work fills this gap timely and provides an attractive alternative to establish optimal sub-models for model averaging.
Subsequently, a $1-\alpha$ confidence interval for the ATE is constructed by the easy-to-implement normal approximation approach.
Numerical evidences shows clearly the superior prediction performance of our proposal over other competing approaches.
We provide additional numerical results, technical lemmas and proofs in the Supplementary Material.


\section{Methodology}\label{sect.2}
\subsection{Notation}\label{sect.2.1}

We consider a panel data setting $\big\{ Y_{it} \big\}_{i=1,t=1}^{N,T}$ with $N$ units and $T$ time periods.
Let $Y_{it}(1)$ and $Y_{it}(0)$ be the outcome of unit $i$ at time $t$ with and without the treatment, respectively. The
treatment effect of the unit $i$ at a post-treatment time-point $t$ is defined as $\Delta_{it}=Y_{it}(1)-Y_{it}(0)$.
However, for any $i$ and $t$, we do not simultaneously observe $Y_{it}(1)$ and $Y_{it}(0)$. Note that one only observes $Y_{it}(1)$ or $Y_{it}(0)$ depending on whether the unit $i$ receives a treatment at time $t$.
The observed outcomes can be usually described as $Y_{it}=d_{it}Y_{it}(1)+(1-d_{it})Y_{it}(0)$, where $d_{it}$ is a binary variable: $d_{it} = 1$ if the unit $i$ receives the treatment at time $t$ and $d_{it} = 0$ otherwise.

Without loss of generality, following the framework of \cite{LB17}, \cite{L20}, \cite{MM21} and \cite{BFR21}, we assume that only the first unit receives a treatment at time $T_1+1$ ($1<T_1<T$) and the remaining $ N-1$ units do not receive a treatment at any time, often referred to as the controls.
Under such an assumption, the observed
outcomes are
\begin{align} \label{eqfoc_A}
Y_{it}=\begin{cases}
Y_{it}(1),&\text{for $i=1$ and $t= T_1+1,\cdots,T$},\\
Y_{it}(0), &\text{for $i=1$ and $t=1,\cdots,T_1$},\\
Y_{it}(0), &\text{for $i=2,\cdots,N$ and $t=1,\cdots,T$.}
\end{cases}
\end{align}
Our main goal is to estimate the ATE parameter of the first unit: $\Delta_{1}=\mathbb{E} \big(\Delta_{1t}\big)= \mathbb{E}\big(Y_{1t}(1)-Y_{1t}(0)\big)$.
The major challenge is that the counterfactual outcome $Y_{1t}(0)$ cannot be observed for the post-treatment time periods $t=T_1+1,\cdots,T$.

Suppose that $\hat \mu_{1t}$ is a generic estimator of the counterfactual outcome $Y_{1t}(0)$ for $t=T_1+1,\cdots,T$. Then one can estimate $\Delta_{1}$ consistently by averaging over the post-treatment period
\begin{align} \label{eqA0}
\hat \Delta_{1}=T_2^{-1}\sum_{t=T_1+1}^T \hat\Delta_{1t} =T_2^{-1}\sum_{t=T_1+1}^T \left(Y_{1t}-\hat \mu_{1t}\right),
\end{align}
where $ T_2 \triangleq T - T_1$ is the post-treatment sample size.

\subsection{Model set-up}\label{sect.2.1A}

The basic approach for constructing the counterfactuals usually depends on the correlations among cross-sectional
units. One often assumes the correlations among cross-sectional units are due to some latent factors.
Following \cite{HCW12,LB17,B21}, we assume that the potential outcome without treatment is generated by
a factor model of the form
\begin{align} \label{eqfoc}
Y_{it}(0)=a_i+\bm b_i^\top\bm f_t+u_{it}, i=1,\cdots,N,t=1,\cdots,T,
\end{align}
where $a_i$ is a subject-specific intercept, $\bm b_i$ is a $\mathbb{K}\times 1$ vector of
factor loadings for the $i$-th unit, $\bm f_t$ is a $\mathbb{K}\times 1$ unobservable common factors
and $u_{it}$ is a weakly dependent and weakly stationary error term.
Write $\bm Y_{t}(0)=(Y_{1t}(0),\cdots, Y_{Nt}(0))^\top$, $ \bm a=(a_1,\cdots, a_N)^\top$, $\bm B=(\bm b_1,\cdots, \bm b_N)^\top$ and
$\bm u_{t}=(u_{1t},\cdots, u_{Nt})^\top$. Then, we can express model \eqref{eqfoc} in the following matrix form
\begin{align} \label{eqfoc*}
\bm Y_{t}(0)=\bm a+\bm B\bm f_t+\bm u_{t}, t=1,\cdots,T,
\end{align}
where the error vector $\bm u_{t}$ satisfies $\mathbb{E}\big(\bm u_{t}\big)=\bm 0$, $\mathbb{E}\big(\bm u_{t}\bm f_t^\top\big)=\bm 0$, $\mathbb{E}\big(u_{it}|d_{1t}\big)=0$ for all $i\neq1$ and $\mathbb{E}\big(\bm u_{t}\bm u_{t}^\top\big)$ is an $N\times N$ diagonal matrix.
These assumptions are similar to those used in \cite{HCW12,LB17}.
Further, all rows of $\bm B$ are bounded (in $L_2$-norm) uniformly by a positive constant and $Rank(\tilde {\bm B})=\mathbb{K}$, where $\tilde {\bm B}$ is the $(N-1) \times \mathbb{K}$ matrix obtained by removing the first row of $\bm B$.
The assumption on $\bm B$ is essential for the calculations leading to \eqref{eqAA}.
\cite{HCW12} proposed to estimate the counterfactual
outcome without estimating unobserved factors
and loadings. Specifically, let $\bm \gamma=(1,-\tilde{\bm\gamma}^\top)^\top$ be an $N\times 1$ vector lying in the null space of $\bm B$, i.e. $\bm \gamma\in \mathcal{N}(\bm B)$. Then, $\bm \gamma^\top\bm B=\bm 0$ \citep{HCW12}. Hence, multiplying \eqref{eqfoc*} by $\bm \gamma^\top$ and re-arranging the terms on two sides leads to
\begin{align} \label{eqAA}
Y_{1t}(0)=\bm \gamma^\top\bm a+ \tilde{\bm\gamma}^\top\tilde{\bm Y}_{t}(0) +\bm \gamma^\top\bm u_{t}
\triangleq\alpha+ \bm X_t^\top\bm \beta+\epsilon_{t}, ~t=1,\cdots,T,
\end{align}
where $\bm X_t \triangleq \big(X_{1t},\cdots,X_{(N-1)t}\big)^\top=\tilde{\bm Y}_{t}(0)=\left(Y_{2t}(0),\cdots, Y_{Nt}(0)\right)^\top$ is an $(N-1)\times 1$ vector composed of the $N-1$ control units' pre-treatment outcomes and $\alpha=\bm \gamma^\top\bm a$, $\bm \beta=\bm \tilde{\gamma}$, $\epsilon_{t}=\bm \gamma^\top\bm u_{t}$.
 We note that equation \eqref{eqAA} holds for
the entire period and $Y_{it}=Y_{it}(0)$ by \eqref{eqfoc_A} for any $i$ and any $t\le T_1$.
Thus, at the pre-treatment period, equation \eqref{eqAA} reduces to
\begin{align} \label{eqA}
Y_{1t}=\alpha+ \bm X_t^\top\bm \beta+\epsilon_{t}\triangleq\mu _{1t}+\epsilon_{t}, ~t=1,\cdots,T_1.
\end{align}

In model (\ref{eqA}), the pre-treatment outcomes of the first unit are regarded as the response, and the pre-treatment outcomes of the $N-1$ control units serve as regressors. Now, we estimate the counterfactual outcome by $\hat\mu_{1t}=\hat\alpha+\bm X_t^\top\hat{\bm\beta}$ for $t=T_1+1,\cdots,T$, where $\hat\alpha$ and $\hat{\bm\beta}$ are generic estimators (i.e., SCM) from fitting model (\ref{eqA}).

\subsection{Robust rank-based model averaging via synthetic control}\label{sect.2.2}
 We provide a model averaging strategy to handle the uncertainty of control units and obtain a robust estimation of the ATE simultaneously.
Suppose that a total of $S$ candidate models $\mathcal{M}_1,\cdots,\mathcal{M}_{S}$ are considered. For $s=1,\cdots,S$, let $\mathcal{A}_s\subseteq \{1,\cdots,N-1\} $ be the index set of covariates to be included in model $\mathcal{M}_s$ with cardinality $|\mathcal{A}_s|=p_s$. Then, in the same spirit as \eqref{eqA}, we consider the $s$th candidate model as
\begin{align} \label{eq0}
\mathcal{M}_s: Y_{1t}=\alpha^{(s)}+ \bm X_t^{(s)\top}\bm\beta^{(s)}+\epsilon_{t}^{(s)},t=1,\cdots,T_1,
\end{align}
where $\bm X_{t}^{(s)}=\left(X_{jt}: j\in \mathcal{A}_s\right)^\top$ is the $p_s$-dimensional covariate vector and
$\alpha^{(s)}$, $\bm\beta^{(s)} =\left(\beta_{1}^{(s)},\cdots, \beta_{p_s}^{(s)}\right)^\top$ are the unknown regression coefficients and $\epsilon_{t}^{(s)}$ is the random error for model $\mathcal{M}_s$.

We first estimate the unknown parameters in the model $\mathcal{M}_s$ for $s=1,\cdots,S$. To safeguard against outliers, we propose to estimate $\bm \beta^{(s)}$ by minimizing the following weighted rank loss
\begin{align} \label{eq1}
\hat{\bm\beta}^{(s)} =&\mathop {\arg \min }\limits_{\bm\beta^{(s)}\in\mathbb{R}^{p_s}} \mathcal{L}\left(\bm\beta^{(s)}\right) \nonumber\\
=&\mathop {\arg \min }\limits_{\bm\beta^{(s)}\in\mathbb{R}^{p_s}} \frac{1}{{{T_1}\left( {{T_1} - 1} \right)}}\sum_{t =1}^{T_1}\sum_{t' \neq t} { b_{tt'}^{(s)}\left| e_t^{(s)}\left(\bm\beta^{(s)}\right) -e_{t'}^{(s)}\left(\bm\beta^{(s)}\right)\right|}\nonumber\\
=&\mathop {\arg \min }\limits_{\bm\beta^{(s)}\in\mathbb{R}^{p_s}} \frac{1}{{{T_1}\left( {{T_1} - 1} \right)}}\sum_{t =1}^{T_1}\sum_{t' \neq t} {b_{tt'}^{(s)}\left| Y_{1tt'}-\bm X_{tt'}^{(s)\top}\bm\beta^{(s)}\right|},
\end{align}
where $e_t^{(s)}\left(\bm\beta^{(s)}\right)=Y_{1t}-\bm X_t^{(s)\top}\bm\beta^{(s)}$, $Y_{1tt'}=Y_{1t}-Y_{1t'}$, $\bm X_{tt'}^{(s)}=\bm X_t^{(s)}-\bm X_{t'}^{(s)}$ and $b_{tt'}^{(s)}$ is a positive and symmetric weight.
Incorporation of the weight $b_{tt'}^{(s)}$ can bound the effect of leverage points in the covariate space. Although many different types of $b_{tt'}^{(s)}$ are recommended in the literature, we will adopt the well-known Mallows-type weight
\begin{align*}
b_{tt'}^{(s)}=h\big(\bm X_t^{(s)}\big)h\big(\bm X_{t'}^{(s)}\big),
\end{align*}
where $h\big(\bm X_t^{(s)}\big)={\rm min}\left\{1,\nu_s/\big[(\bm X_t^{(s)}-\hat{\bm\varsigma}^{(s)})^\top \big(\hat{\bm \Xi}^{(s)}\big)^{-1} (\bm X_t^{(s)}-\hat{\bm\varsigma}^{(s)})\big]\right\}$, $\nu_s$ is the $95\%$ quantile of the chi-square distribution with $p_s$ degrees of freedom, and $\hat{\bm\varsigma}^{(s)}$ and $\hat{\bm \Xi}^{(s)}$ are the robust minimum volume ellipsoid estimators
of $\mathbb{E}\big(\bm X_t^{(s)}\big)$ and ${\rm Cov}\big(\bm X_t^{(s)}\big)$ respectively \citep{Rvan90,HM98,WL09}. Such weights can be easily implemented by the function ``cov.rob" in R package `` {\sf MASS}".

Furthermore, we notice that the objective function $\mathcal{L}\left(\bm\beta^{(s)}\right) $ does not depend on the intercept since $\alpha^{(s)}$ is canceled out. Following \cite{WL09,JLLC24}, we estimate $\alpha^{(s)}$ by the sample median approach as $\hat\alpha^{(s)}=\mbox{Median} \left\{e_t^{(s)}\big(\hat{\bm\beta}^{(s)}\big): 1 \leq t \leq T_1  \right\}$.

Then, under the $s$th candidate model $\mathcal{M}_s$, the predicted counterfactual outcome of $Y_{1t}(0)$ is given by
\begin{align} \label{eq7A}
\hat \mu_{1t}^{{(s)}}=\hat\alpha^{(s)}+\bm X_t^{(s)\top}\hat{\bm\beta}^{(s)} ~ {\rm for} ~t=T_1+1,\cdots,T.
\end{align}
When the estimators for counterfactual outcomes $Y_{1t}(0)$ under all candidate models are readily available, the model averaging estimator for the ATE is given by
\begin{align} \label{eq7}
\hat\Delta_{1}(\bm w)=T_2^{-1}\sum_{t=T_1+1}^T \hat\Delta_{1t}(\bm w)=T_2^{-1}\sum_{t=T_1+1}^T \left(Y_{1t}- \hat\mu_{1t}(\bm w) \right),
\end{align}
where $\hat\mu_{1t}(\bm w)=\sum_{s=1}^{S} {w_s\hat \mu_{1t}^{{(s)}} }$ is the estimating the counterfactual outcome and $\bm w=\big(w_1,\cdots,w_S\big)^\top\in \mathcal{W}$ is the model weight vector with $\mathcal{W}=\big\{ \bm w\in [0,1]^{S}:  \sum\limits_{s=1}^{S} {w_s}=1 \big\}$.

\begin{remark}\label{rk1}
An appealing feature of minimizing $\mathcal{L}\left(\bm\beta^{(s)}\right) $ is that its computation can be conveniently carried out using the standard least absolute deviation regression. Let us consider an augmented data set  $\big\{ \big(\tilde Y_{1t^\star}, \tilde{\bm X}_{t^\star}^{(s)} \big):1\leq t^\star\leq T_1(T_1-1) \big\} =\big\{\big( b_{tt'}^{(s)} Y_{1tt'},b_{tt'}^{(s)}\bm X_{tt'}^{(s)}\big):1\leq t'\neq t\leq T_1\big\}$. Then, we can rewrite $\mathcal{L}\left(\bm\beta^{(s)}\right)$ as $\sum_{t^\star =1}^{T_1(T_1-1)} \left|\tilde Y_{1t^\star}-  \tilde{\bm X}_{t^\star}^{(s)\top}\bm\beta^{(s)} \right|$, and thus $\hat{\bm\beta}^{(s)}$ can be implemented easily by the function {\sf rq} in R package {\sf quantreg}.
\end{remark}

\subsection{Estimation of optimal model weights}\label{sect.2.3}

The weights in (\ref{eq7}) play a central role in producing reliable estimation results.
Although the $K$-fold cross-validation has been widely used to determine model weights (e.g., \cite{GLW15,LS15,ZL23}), it may be problematic for time series or dependent data.
The traditional $K$-fold cross-validation randomly shuffles and partitions observations, which breaks the inherent temporal dependence in time series and introduces non-negligible bias \citep{hyndman2018forecasting}. To avoid the issue, we adopt the growing window cross-validation (GWCV) to select the model weights. This criterion strictly preserves chronological order, ensuring that model training is always based on past observations to predict future outcomes. Specifically, we set $\mathcal{I}_t^c=\{1,2,\cdots,t-1\}$ and $\mathcal{I}_t=\{t\}$ for $t=\frac{T_1}{2}+1,\frac{T_1}{2}+2,\cdots,T_1$.
Let $\alpha^{(s)}_{[t]}$ and $\hat{\bm\beta}^{(s)}_{[t]}$ be the estimators of $\alpha^{(s)}$ and $\bm\beta^{(s)}$ in the $s$-th candidate model by using the observations from the training set $\mathcal{I}_t^c$, that is
\begin{align} \label{eq0p}
\hat{\bm\beta}^{(s)}_{[t]} =\mathop {\arg \min }\limits_{\bm\beta^{(s)}\in\mathbb{R}^{p_s}} \mathcal{L}_{[t]}\left(\bm\beta^{(s)}\right) ~{\rm and}~ \hat\alpha_{[t]}^{(s)}=\mbox{Median} \left\{e_l^{(s)}\big(\hat{\bm\beta}^{(s)}_{[t]}\big): l\in \mathcal{I}_t^c  \right\},
\end{align}
where $\mathcal{L}_{[t]}\left(\bm\beta^{(s)}\right)=\frac{1}{{{|\mathcal{I}_t^c|}\left( {{|\mathcal{I}_t^c|} - 1} \right)}}\sum\limits_{l \in \mathcal{I}_t^c}\sum\limits_{l' \neq l\in \mathcal{I}_t^c } {b_{ll'}^{(s)}\left| Y_{1ll'}-\bm X_{ll'}^{(s)\top}\bm\beta^{(s)}\right|}$.

For the $s$-th candidate model, the prediction for the observation in the test set $\mathcal{I}_t$ is given by: $\hat \mu_{1[t]}^{(s)}=\hat\alpha_{[t]}^{(s)}+\bm X_t^{(s)\top}\hat{\bm\beta}_{[t]}^{(s)} $, $t=\frac{T_1}{2}+1,\cdots,T_1, s=1,\cdots,S$.
Write $\bm{\hat \mu}_{1[t]}=\left(\hat \mu_{1[t]}^{(1)},\cdots,\hat \mu_{1[t]}^{(S)}\right)^\top$ and $\bm{\hat \mu}_{1} =\left(\bm{\hat \mu}_{1[t]}:t=\frac{T_1}{2}+1,\frac{T_1}{2}+2,\cdots,T_1\right)^\top$.
To reduce the influence of outliers,
the Mallows-type weights $\tau_{tt'}$ are adopted again, given by
\begin{align}\label{MMW}
\tau_{tt'}=\hbar\big(\bm{\hat \mu}_{1[t]}\big)\hbar\big(\bm{\hat \mu}_{1[t']}\big),
\end{align}
where $\hbar\big(\bm{\hat \mu}_{1[t]}\big) ={\rm min}\left\{1,\upsilon/\big[(\bm{\hat \mu}_{1[t]}-\hat{\bm\zeta})^\top \hat{\bm \Sigma}^{-1} (\bm{\hat \mu}_{1[t]}-\hat{\bm\zeta})\big]\right\}$, and $\upsilon$ is the $95\%$ quantile of the chi-square distribution with $S$ degrees of freedom. Here $\hat{\bm\zeta}$ and $\hat{\bm \Sigma}$ are the robust minimum volume ellipsoid estimators of $\mathbb{E}\big(\bm{\hat \mu}_{1[t]}\big)$ and ${\rm Cov}\big(\bm{\hat \mu}_{1[t]}\big)$ respectively.

We attain $\hat {\bm w}=\mathop {\arg \min }\limits_{\bm w \in\mathcal{W}}\mathcal{C}\left(\bm w\right)=\mathop {\arg \min }\limits_{\bm w \in\mathcal{W}}\frac{4}{T_1(T_1-2)}\sum\limits_{t =\frac{T_1}{2}+1}^{T_1}\sum\limits_{t'\neq t}\tau_{tt'}\left| Y_{1tt'}-\bm{\hat \mu}_{1[tt']}^\top\bm w\right|$, where $\bm{\hat \mu}_{1[tt']}=\bm{\hat \mu}_{1[t]}-\bm{\hat \mu}_{1[t']}$. Using similar arguments of Remark \ref{rk1}, the loss function $\mathcal{C}\left(\bm w\right)$ can be rewritten as $\mathcal{C}\left(\bm w\right)=\frac{4}{T_1(T_1-2)}\sum\limits_{t^\star =1}^{T_1(T_1-2)/4} \left|\bar Y_{1t^\star}-  \bm{\bar \mu}_{1[t^\star]}^\top \bm w\right|$ by using an augmented data set $\big\{ \big(\bar Y_{1t^\star}, \bm{\bar \mu}_{1[t^\star]} \big):1\leq t^\star\leq \frac{T_1}{2}(\frac{T_1}{2}-1) \big\} =\big\{\big( \tau_{tt'} Y_{1tt'},\tau_{tt'}\bm{\hat \mu}_{1[tt']}\big):\frac{T_1}{2}+1\leq t'\neq t\leq T_1\big\}$.
Although a general nonlinear constrained solver (e.g., the {\sf donlp} in R package {\sf Rsolnp}) can be used to obtain the minimizer of $\mathcal{C}\left(\bm w\right)$, it maybe time-consuming and prone to issues regarding local minima and solver dependency.
Note that the objective function $\mathcal{C}\left(\bm w\right)$ is convex in $\bm w$.
Hence, we can convert the constrained $L_1$-loss minimization into an equivalent linear program by introducing non-negative slack variables. We can obtain the optimal weight vector $\hat{\bm{w}}$ by solving the following standard linear program:
\[
\begin{aligned}
&\min_{\bm{w},\bm{e}}  \quad \frac{4}{T_1(T_1-2)} \sum_{t^\star=1}^{T_1(T_1-2)/4} e_{t^\star} \\
\text{s.t.} \quad
&\begin{cases}
\bar{Y}_{1t^\star}-\bar{\bm{\mu}}_{1[t^\star]}^\top \bm{w} \leq e_{t^\star}, \\
-\left(\bar{Y}_{1t^\star}-\bar{\bm{\mu}}_{1[t^\star]}^\top \bm{w} \right)\leq e_{t^\star}, \\
e_{t^\star} \ge 0,
\sum_{s=1}^{S} w_s = 1, 
w_s \ge 0,
\end{cases}
\end{aligned}
\]
where $\bm w=\big(w_1,\cdots,w_S\big)^\top$ and $\bm e=\left(e_{1},\cdots,e_{T_1(T_1-2)/4}\right)^\top$.

Plugging $\hat {\bm w}$  into (\ref{eq7}), we obtain the robust rank-based model averaging ({\sf RRMA}) estimator for the ATE as $\hat\Delta_{1}(\hat {\bm w})$.

\section{Theoretical results}\label{sect.3}
\subsection{Asymptotic properties when all candidate models are mis-specified}

In this section, we assume all candidate models are mis-specified. Let $\lambda _{\min }(\bm A)$ and $\lambda _{\max }(\bm A)$ be the minimum and maximum eigenvalues of a symmetric matrix $\bm A$. For a vector $\bm \xi=(\xi_1,\cdots,\xi_s)^\top\in {\mathbb{R}^s}$, denote its $L_2$ norm by ${\left\| \bm \xi \right\|_2} = {\left( {\xi_1^2 +  \cdots  + \xi_s^2} \right)^{1/2}}$.
Furthermore, we define the population pseudo-true parameters of $\alpha ^{(s)}$ and $\bm \beta ^{(s)}$ under the model $\mathcal{M}_s$ by
${\bm\beta ^{*(s)}} =  \mathop {\arg \min }\limits_{{\bm\beta ^{(s)}} \in {\mathbb{R}^{{p_s}}}}\mathbb{E}\left[ \frac{1}{{T_1(T_1 - 1)}}{\sum\limits_{t = 1}^{{T_1}} {\sum\limits_{t' \ne t} {b_{tt'}^{(s)}\Big| {{Y_{1tt'}} - \bm X_{tt'}^{(s) \top }{\bm\beta ^{(s)}}} \Big|} } } \right] $
and ${\alpha ^{*(s)}} = \mathop {\arg \min }\limits_{\alpha ^{(s)}\in\mathbb{R}} \mathbb{E}\left[\frac{1}{{T_1}} {\sum\limits_{t = 1}^{{T_1}} {\left| { {{Y_{1t}} - \bm X_t^{(s) \top }{\bm\beta ^{*(s)}}}  - {\alpha ^{(s)}}} \right|} } \right].$
Let $G_{\epsilon_{t}}\left(\cdot|\bm X_t\right)$ and $F_{\epsilon_{tt'}}\left(\cdot|\bm X_{tt'}\right)$ be the respective conditional cumulative distribution functions of $\epsilon_{t}$ given $\bm X_{t}$ and $\epsilon_{tt'}$ given $\bm X_{tt'}$ and $g_{\epsilon_{t}}\left(\cdot|\bm X_t\right)$ and $f_{\epsilon_{tt'}}\left(\cdot|\bm X_{tt'}\right)$ be the corresponding conditional density functions, where $\epsilon_{tt'}\triangleq \epsilon_t-\epsilon_{t'}$.
Let $\bm\Omega ^{(s)} =  \mathop {\lim }\limits_{{T_1} \to \infty } \mathbb{E}\Big\{\frac{1}{{{T_1}\left( {{T_1} - 1} \right)}}\sum\limits_{t = 1}^{{T_1}} {\sum\limits_{t' \ne t} {b_{tt'}^{(s)}{f_{{\epsilon_{tt'}}}}} } \big( -u_{tt'}^{*(s)} | {\bm X_{tt'}}   \big)\bm X_{tt'}^{(s)}\bm X_{tt'}^{(s) \top } \Big\}$, $\bm\Psi ^{(s)} =\mathop {\lim }\limits_{{T_1} \to \infty } {\rm Cov}\Big( \frac{\sqrt{T_1}}{ {{T_1}\left( {{T_1} - 1} \right)} }\sum\limits_{t = 1}^{{T_1}} {\sum\limits_{t' \ne t} {b_{tt'}^{(s)} \psi \left( {{Y_{1tt'}},\bm X_{tt'}^{(s) \top }{\bm \beta ^{*(s)}}} \right) \bm X_{tt'}^{(s)}} }  \Big)$, $u_{tt'}^{*(s)} = u_t^{*(s)} - u_{t'}^{*(s)} $ and $u_t^{*(s)} = {{ \alpha }} + \bm X_t^ \top {{\bm\beta}} - {\alpha ^{*(s)}} - \bm X_t^{(s) \top }{\bm\beta ^{*(s)}}$ and
$\psi \left( y,x \right) \triangleq \frac{\partial \mid y-x\mid}{\partial x}= -1 + 2{\rm I}\big( y-x \leq 0 \big)$ with ${\rm I}\left(\cdot\right)$ being an indicator function. Moreover, let
$\mathcal{S}_1^{(s)}=\frac{\sqrt{T_1}}{T_1}\sum\limits_{t = 1}^{{T_1}} {\left[ {\mathrm{I}\left( {{\epsilon_t} + u_t^{*(s)} \leq 0} \right) - 1/2} \right]}$, $\bm {\mathcal{S}}_2^{(s)}=\frac{\sqrt{T_1}}{{{T_1}\left( {{T_1} - 1} \right)}}\sum\limits_{t = 1}^{{T_1}} {\sum\limits_{t' \ne t} {b_{tt'}^{(s)}\left[ {\mathrm{I}\left( {{\epsilon_{tt'}} + u_{tt'}^{*(s)} \leq 0} \right) - 1/2} \right]} }\bm  X_{tt'}^{(s)}$, $\bm {\mathcal{S}}^{(s)}= \left(\mathcal{S}_1^{{(s)}},\bm {\mathcal{S}}_2^{(s)\top} \right)^\top$ and
$\bm {\mathcal{S}}=\left(\bm {\mathcal{S}}^{(1)\top},\cdots,\bm {\mathcal{S}}^{(S)\top} \right)^\top$.
Define $\bm\theta^{*(s)}=\left(\alpha^{*(s)} ,\bm\beta^{*(s)\top}\right)^\top$, $\bm\theta^*=\left(\bm\theta^{*(1)\top},\cdots,\bm\theta^{*(S)\top}\right)^\top$, $\hat{\bm\theta}^{(s)}=\left(\hat\alpha^{(s)} ,\hat{\bm\beta}^{(s)\top}\right)^\top$ and $\hat{\bm\theta}=\left(\hat{\bm\theta}^{(1)\top},\cdots,\hat{\bm\theta}^{(S)\top}\right)^\top$.
Let $\bar p = \mathop {\max }\limits_{1 \le s \le S} {p_s}$ be the maximum number of parameters in all candidate models.
Meanwhile, following \cite{LB17,L20}, we assume that $\phi  = \mathop {\lim }\limits_{{T_1},{T_2} \to \infty } \sqrt {{T_2}/{T_1}} $ is a finite nonnegative constant. We also allow $\bar p$ and $S$ to diverge to infinity as $T_1\rightarrow \infty$.
We impose the following technical conditions:

\begin{enumerate}
\item[(C1)] $F_{\epsilon_{tt'}}(\cdot\mid \boldsymbol{X}_{tt'})$ and $G_{\epsilon_{t}}(\cdot\mid \boldsymbol{X}_{t})$ are absolutely continuous. Their corresponding probability density functions $f_{\epsilon_{tt'}}(\cdot\mid \boldsymbol{X}_{tt'})$, $g_{\epsilon_{t}}(\cdot\mid \boldsymbol{X}_{t})$, as well as their first-order derivatives $\dot{f}_{\epsilon_{tt'}}(\cdot\mid \boldsymbol{X}_{tt'})$ and $\dot{g}_{\epsilon_{t}}(\cdot\mid \boldsymbol{X}_{t})$, are continuous and uniformly bounded.

\item[(C2)]  Assume that there exist two positive constants $\rho_{1}^{(s)}$
and $\rho_{2}^{(s)}$ such that ${\lambda _{\min }}\left( {{\bm\Omega ^{(s)}}} \right)\geq \rho_{1}^{(s)}$ and ${\lambda _{\max }}\left( {{\bm\Psi^{(s)}}} \right)\leq \rho_{2}^{(s)}$ for $s=1,\cdots, S$.

\item[(C3)] There exists a positive constant $c$ such that $| X_{jt}|<c$ for $1\leq t\leq T_1$ and $1\leq j\leq N-1$.

\item[(C4)] Assume $\rho_1^{-2}\rho_2=o(T_1\log T_1)$, $\bar p^4\log T_1 =o(\rho_1^2 T_1)$ and
   $ST_1^{- \rho_1^{2}\rho_2^{-1}\bar p }=o(1)$, where $\rho_1=\mathop {\min }\limits_{1 \le s \le S} \rho_{1}^{(s)}$ and $\rho_2=
\mathop {\max }\limits_{1 \le s \le S} \rho_{2}^{(s)}$.

\item[(C5)] Let $\varrho_t=(Y_{1t},Y_{2t},\ldots,Y_{Nt},\Delta_{1t}d_{1t})$ for $t = 1,\ldots,T$, where $d_{1t} = 0$ if $t\leq T_1$and $d_{1t} = 1$ if $t\geq T_1 + 1$. Assume that $\{\varrho_t\}_{t = 1}^{T_1}$ and $\{\varrho_t\}_{t = T_1+1}^{T}$ are both weakly dependent stationary processes. Define $\rho(\tau)=\max\limits_{1\leq t\leq T}\max\limits_{1\leq i,j\leq N + 1}|{\rm Cov}(\varrho_{it},\varrho_{j,t+\tau})/\sqrt{\text{var}(\varrho_{it})\text{var}(\varrho_{j,t+\tau})}|$. Then there exists some finite positive constants $C>0$ and $0<\lambda<1$ such that $\rho(\tau)\leq C\lambda^{\tau}$.

\item[(C6)] For any $ \bm a\in \mathbb{R}^{\sum_{s=1}^S p_s+S}$ such that $\|\bm a\|_2=1$,
$\bm a^\top \bm\Xi ^{-1/2} \bm {\mathcal{S}}  \mathop  \to\limits^d \mathcal{N}\left(0,1 \right)$, where $\bm \Xi=
\mathop {\lim }\limits_{{T_1} \to \infty } \mathbb{E} \{\bm {\mathcal{S}} \bm {\mathcal{S}}^\top \}$ is a positive definite matrix.

\end{enumerate}

Condition (C1) assumes that $f_{\epsilon_{tt'}}\left(\cdot|\bm X_{tt'}\right)$ and $\dot f_{\epsilon_{tt'}}\left(\cdot|\bm X_{tt'}\right)$, $ g_{\epsilon_{t}}\left(\cdot|\bm X_t\right)$ and $\dot g_{\epsilon_{t}}\left(\cdot|\bm X_t\right)$ are uniformly bounded away from zero and infinity, which are commonly imposed in robust regression \citep{LS15,JLLC24}. Condition (C2) imposes some restrictions on the minimum and maximum eigenvalues of ${{\bm\Omega ^{(s)}}}$ and ${{\bm\Psi^{(s)}}}$ respectively, and are often used to establish the upper bound of estimation error. Condition (C3) assumes that the predictor $X_{jt}$ has bounded absolute
value, usually considered in robust rank-based regression analysis, see Condition C1 of \cite{HMZZ22} and Assumption 6 of \cite{FHXZZ}.
Condition (C4) gives some mild restrictions on $S$, $T_1$ and $\bar p$, which is analogous to Assumption A.3 of \cite{LS15}.
Condition (C5) assumes that the data is a weakly stationary $\rho$-mixing process with an exponential decay
rate. Many weakly dependent processes, including certain strictly stationary ARMA processes, are recognized as $\rho$-mixing processes that exhibit exponential decay rates. Moreover, this condition ensures that the estimator obtained from the pretreatment samples is asymptotically
independent with a quantity that involves the post-treatment samples, which is similar to Assumption 4 of \cite{L20}.
From the definition of ${\bm\beta ^{*(s)}}$ and ${\alpha ^{*(s)}}$, it is easy to verify $\mathbb{E}\left(\bm {\mathcal{S}}^{(s)} \right)=\bm 0 $. Moreover, $\bm {\mathcal{S}}$ are expressed as the sum of a series of weakly
dependent variables. Thus, Condition (C6) is rational to assume that
$\bm a^\top \bm\Xi ^{-1/2} \bm {\mathcal{S}}$ converges to a standard normal distribution by the central limit theorems, because laws of large numbers
and central limit theorems are known to hold for many weakly
dependent and weakly stationary processes and $p_s$ diverges slower than $O(T_1^{1/4})$ by Condition (C4).
Similar assumptions can refer to Assumptions 1 and 2 of \cite{L20}, Assumption 4 of \cite{LB17}.

\begin{proposition}\label{pro.1}
Under conditions (C1)--(C4) hold, it holds that
$\mathop {\max }\limits_{1 \le s \le S} {\left\| {{\bm{\hat \beta }^{(s)}} - {\bm\beta ^{*(s)}}} \right\|_2} = O_p\left( {\sqrt {\bar p\log T_1/{T_1}} } \right)$.
\end{proposition}

\begin{proposition}\label{pro.2}
Suppose conditions (C1)--(C6) hold. Then, for any $ \bm a\in \mathbb{R}^{\sum_{s=1}^S p_s+S}$ such that $\|\bm a\|_2=1$, we have
\begin{align*}
  \sqrt {{T_1}} \bm a^\top \left(\bm\Phi \bm\Xi \bm\Phi^\top\right)^{-1/2} \left( \hat{\bm\theta}-{\bm\theta}^* \right)   \mathop  \to \limits^d
\mathcal{N}\left( 0, 1\right),
 \end{align*}
 where $\bm\Phi={\rm diag}\left(\bm\Phi^{(1)},\cdots,\bm\Phi^{(S)}\right)$,
 $\bm\Phi^{(s)}=\left( {\begin{array}{*{20}{c}}
  \phi_{11}^{(s)} &  \bm\phi_{12}^{(s)}  \\
  \bm\phi_{21}^{(s)} &  \bm\phi_{22}^{(s)}  \\
\end{array}} \right)$,
 $\phi_{11}^{(s)}= \frac{-1}{{\mathbb{E}\left[ {g_{\epsilon_{t}}\left( { - u_t^{*(s)}}|\bm X_t \right)} \right]}} $, $\bm\phi_{12}^{(s)}=\frac{{\mathbb{E}\left[{g_{\epsilon_{t}}\left( { - u_t^{*(s)}}|\bm X_t \right)\bm X_t^{(s)\top}} \right]\left(  \bm\Omega ^{(s)} \right)^{ - 1} } }{2\mathbb{E}\left[ {g_{\epsilon_{t}}\left( { - u_t^{*(s)}} |\bm X_t\right)} \right]}$, $\bm \phi_{21}^{(s)}=\bm 0$ and $\bm\phi_{22}^{(s)}=-\left( {{\bm \Omega^{(s)}}} \right)^{ - 1}$.
\end{proposition}

\begin{remark}\label{rk2}
Existing theoretical results (such as consistency and asymptotic distribution) on rank-based estimation usually hold under a correctly specified model \citep{WL09,WPBLW20}. However, a common assumption in the field of model averaging is that all candidate models are wrong. Therefore, existing theoretical results cannot be directly utilized here. When all the candidate models $\{\mathcal{M}_s, s=1,\cdots, S\}$ are potentially mis-specified, Proposition \ref{pro.1} ensures the uniform convergence of $\bm{\hat \beta }^{(s)}$, and Proposition \ref{pro.2} provides the limiting distribution of submodel estimates $ \hat{\bm\theta}$.
Consequently, our theoretical framework does not need the assumption of correctly specified models. This confers greater robustness and generality, extendable to a wide range of real-world scenarios where model mis-specification is inevitable.
\end{remark}

Denote
$\mathcal{R}\left( \bm w \right) = \mathbb{E}\Big\{ \frac{4}{T_1(T_1-2)}\sum\limits_{t =\frac{T_1}{2}+1}^{T_1}\sum\limits_{t'\neq t}\tau_{tt'}\left| Y_{1tt'}-\bm{\hat \mu}_{1[tt']}^\top\bm w\right|\Big\}$ as the rank-based risk function.
Further, let $\mathcal{R}^*\left( \bm w \right) $ be the rank-based risk function calculated based on the population pseudo-true parameters instead of the estimated parameters, that is,
$\mathcal{R}^*\left( \bm w \right) = \mathbb{E}\Big\{ \frac{4}{T_1(T_1-2)}\sum\limits_{t =\frac{T_1}{2}+1}^{T_1}\sum\limits_{t'\neq t}\tau_{tt'}\left| Y_{1tt'}-\bm{\mu}_{1tt'}^{*\top}\bm w\right|\Big\} $, where $\bm{\mu}_{1tt'}^*=\bm{\mu}_{1t}^*-\bm{\mu}_{1t'}^*$, $\bm{\mu}_{1t}^*=\left(\mu_{1t}^{*(1)},\cdots,\mu_{1t}^{*(S)}\right)^\top$ and $\mu_{1t}^{*(s)}=\alpha^{*(s)}+\bm X_t^{(s)\top}\bm\beta^{*(s)} $ for $t=\frac{T_1}{2}+1,\cdots,T_1$ and $s=1,\cdots,S$.
The unattainable optimal weight vector is $\bm w^*=\mbox{argmin}_{\bm w\in{\mathcal{W}}}{\mathcal{R}^*(\bm w)}$. We assume that $\bm w^*$ is an interior point of ${\mathcal{W}}$. Denote ${\xi _{{T_1}}} = \mathop {\inf }\limits_{\bm w \in \mathcal{W}} {\mathcal{R}^*}\left(\bm w \right)={\mathcal{R}^*}\left(\bm w^* \right) $. To establish the asymptotic optimality, consistency of $\hat{\bm w}$ and the asymptotic distribution of $\hat\Delta_{1}(\hat {\bm w})$, we further impose the following conditions.

\begin{enumerate}
\item[(C7)] Assume that $\xi _{{T_1}}^{ - 1}\mathop {\sup }\limits_{\bm w \in \mathcal{W}} \frac{4}{T_1(T_1-2)}\sum\limits_{t =\frac{T_1}{2}+1}^{T_1}\sum\limits_{t'\neq t} \tau _{tt'}  \left|\big|    {{Y_{1tt'}} - \bm{\hat \mu} _{1[tt']}^\top \bm w} \big| - \big| {{Y_{1tt'}} - \bm\mu _{1tt'}^{*\top}\bm w} \big| \right|  $ is uniformly integrable.

\item[(C8)] ${\xi _{T_1}^{ - 1}\sqrt{ {\bar p^2}\log T_1/T_1 }}=o(1)$.

\item[(C9)] There exist two positive constants $ c_1$ and $ c_2$ such that
${\rm Pr}\{0<c_1  \leq \lambda _{\min }\left\{ \bm\Upsilon \right\}\leq \lambda _{\max }\left\{ \bm\Upsilon \right\} \leq c_2 S \}$ tends to 1,
where $\bm\Upsilon=  \frac{4}{T_1(T_1-2)}\sum\limits_{t =\frac{T_1}{2}+1}^{T_1}\sum\limits_{t'\neq t} {{\tau _{tt'}}\bm{\mu} _{1tt'}^* \bm{\mu} _{1tt'}^{* \top }} $.

\item[(C10)]
For the post-treatment periods $t=T_1+1,\cdots, T$, $\epsilon_{t}$ and $\Delta_{1t} $ are weakly dependent and weakly stationary processes, and $v_t$ has mean zero and satisfies the central limit theorem: $T_2^{-1/2}\sum_{t=T_1+1}^T v_t\mathop  \to \limits^d \mathcal{Z} \sim \mathcal{N}\left(0,\Sigma_v\right)$, where $v_t= \epsilon_{t}+\Delta_{1t} -\Delta_{1}$ and $\Sigma_v=
\mathop {\lim }\limits_{{T_2} \to \infty } T_2^{ - 1}\sum\limits_{t = T_1+1}^{T}\sum\limits_{s = T_1+1}^{T} \mathbb{E}\left({v_tv_s}\right)$.
\end{enumerate}

Condition (C7) concerns the integrability of some functions and is widely applied to establish the asymptotic optimality of the weight estimator, see Assumption 3 of \cite{FHXZZ}.
Condition (C8) allows that $\bar p$ increases at a rate that
does not exceed $\xi _{T_1}\left(T_1/\log T_1 \right)^{1/2}$. Similar assumptions can refer to Assumption 5 in \cite{ZL23} and Assumption 5 in \cite{FHXZZ}.
Condition (C9) requires that the minimum eigenvalue of $\bm\Upsilon $ is asymptotically bounded and the maximum eigenvalue of $\bm\Upsilon $ is the same order as $S$ (the number of candidate models), which
 is analogous to Condition 3 of \cite{LLWL22}.
Condition (C10) requires that a central limit theorem applies to a partial sum of $v_t$, which is essential to establish the asymptotic distribution of $\hat\Delta_{1}(\hat {\bm w})$. Similar assumptions can refer to Assumption 5 of \cite{LB17}, Assumption 3 of \cite{L20} and Condition (A3) of \cite{RGG24}.

\begin{thm}\label{th.1}
Suppose conditions (C1)--(C4) and (C7)--(C8) hold. Then, it holds that
$\frac{{\mathcal{R}\left( \bm{\hat w} \right)}}{{\mathop {\inf }\limits_{\bm w \in \mathcal{W}} \mathcal{R}\left(\bm w \right)}}\mathop  \to \limits^p 1$.
\end{thm}

\begin{thm}\label{th.2}
 If conditions (C1)--(C4) and (C9) are satisfied, together with $S \bar p^4\log T_1/T_1 \rightarrow 0$ as $T_1\rightarrow \infty$, then there exists a local minimizer $\bm{\hat w}$ of
 $\mathcal{C}\left( \bm w \right)$ such that
${\left\| {\bm{\hat w} - \bm w^*} \right\|_2}=O_p\big(\sqrt {S\bar p^4\log T_1/T_1}  \big).$
\end{thm}

\begin{thm}\label{th.3}
Under the same conditions in Theorem \ref{th.2}, together with conditions (C5), (C6) and (C10), we have
\begin{align*}
 \sqrt {{T_2}} \left( {\hat \Delta _1(\hat{\bm w})- {\Delta _1}} -\sum\limits_{s = 1}^S { w}_s^*\mathbb{E} \left( { u_t^{*(s)}}\right)  \right)
 \mathop  \to \limits^d  \mathcal{N}\left(0,\phi^2\bm{\mathcal{H}}\bm\Phi \bm\Xi \bm\Phi^\top  \bm{\mathcal{H}}^\top +\Sigma_v\right),
\end{align*}
where $\bm{\mathcal{H}}=\left( \bm{\mathcal{H}}^{(1)},\cdots, \bm{\mathcal{H}}^{(S)}\right)$
with $\bm{\mathcal{H}}^{(s)}= \left(w_s^*, w_s^*\mathbb{E}\left(\bm X_t^{(s)}\right)^\top\right)$, $ u_t^{*(s)}=\alpha+ \bm X_t^\top\bm\beta - {\alpha ^{*(s)}} -\bm X_t^{(s)\top}{\bm\beta ^{*(s)}}$ is the approximation bias for the $s$-th candidate model and $\Sigma_v$ is specified in Condition (C10).
\end{thm}

\begin{remark}\label{rk5}
Theorem \ref{th.1} establishes that the weight estimator is optimal in the sense of making the rank-based risk function $\mathcal{R}\left( \bm w \right)$ as small as possible among all feasible weight vectors. Theorem \ref{th.2} indicates that the weight estimator obtained from the GWCV approaches the optimal weight $\bm w^*$ at the rate of $\sqrt {S\bar p^4\log T_1/T_1}$.
Theorem \ref{th.3} demonstrates that $\hat\Delta_1(\hat {\bm w})$ asymptotically follows a normal distribution as $T_1$ and $T_2$ both tend to infinity. It is not surprising that $\hat\Delta_1(\hat {\bm w})$ has a non-ignorable bias, since all candidate models are mis-specified in this section. Furthermore,
Theorem \ref{th.3} extends the classical time-dependent mean model averaging framework of \citet{li2015flexible} to causal inference via a robust rank-based model averaging, offering enhanced robustness to outliers and heavy tails while reducing risk of model mis-specification.
\end{remark}

\subsection{Asymptotic properties when at least one quasi-correct model is in the model set}\label{sect.3.2}

In general not all regressors in \eqref{eqA} contribute non-zero effects. Without loss of generality, let $\mathcal{A}_{true}\triangleq \{j:\beta_j\neq 0,1\leq j \leq N-1\}$ be the true active index set. Recall that $\mathcal{A}_s $ is the index set of covariates to be included in model $\mathcal{M}_s$. The $s$th candidate model $\mathcal{M}_s$ is a \emph{quasi-correct} model if $\mathcal{A}_{true} \subseteq\mathcal{A}_s $. Particularly, when $\mathcal{A}_{true}=\mathcal{A}_s $, the model $\mathcal{M}_s$ is exactly a \emph{correct} model. Actually, the set of quasi-correct models consists of the correct and over-fitted models.
In this section, we consider the scenario where the candidate model set contains at least one quasi-correct model.
Denote $\mathcal{I}_{\rm{cor}}\subseteq \left\{1,2,\cdots, S\right\}$ to be the index set of the
quasi-correct models in the candidate model set with a finite model size.
Further, denote $ w_{\rm{cor}}=\sum\nolimits_{s\in \mathcal{I}_{\rm{cor}}} {w_s}$ and $ \hat w_{\rm{cor}}=\sum\nolimits_{s\in \mathcal{I}_{\rm{cor}}} {\hat w_s}$.
To prove the properties of $\hat{\bm w}$ and $\hat\Delta_{1}(\hat {\bm w})$, an additional condition is needed.
\begin{enumerate}
\item[(C11)] The random error $\epsilon_{t}$ is symmetrically distributed, and the medians of the pairwise difference $\epsilon_{tt'}=\epsilon_{t}-\epsilon_{t'}$ given $\bm X_{tt'}$ are 0 for any $t\neq t'$.
\end{enumerate}

Condition (C11) indicates that  the approximation bias
$u_t^{*(s)} = \alpha  + \bm X_t^ \top {{\bm\beta}}- {\alpha ^{*(s)}} - \bm X_t^{(s) \top }{\bm\beta ^{*(s)}}$ is zero when $\mathcal{M}_s$ is the quasi-correct model. It is an essential assumption to ensure that the proposed estimator $\hat\Delta_{1}(\hat {\bm w})$ for the ATE is consistent, when there is at least one quasi-correct model in
the candidate model set. A similar assumption can refer to \cite{BW05,FHXZZ}.

\begin{thm}\label{th.4}
Suppose the conditions of Theorem \ref{th.2} and condition (C11) hold. If there is at least one quasi-correct model in
the candidate model set, we have $\hat w_{\rm{cor}}\mathop  \to \limits^p 1$.
\end{thm}

\begin{thm}\label{th.5}
Suppose the conditions of Theorem \ref{th.4} and condition (C10) hold. If there is at least one quasi-correct model in
the candidate model set, it holds that $\left|\hat\Delta_{1}(\hat {\bm w})-\Delta_{1}\right|^2={O_p}\left( S^2\bar p^4\log T_1/T_1 +T_2^{-1}\right)$.
\end{thm}

\begin{thm}\label{th.6*}
Suppose the conditions of Theorem \ref{th.3} and condition (C11) hold. If there is at least one quasi-correct model in
the candidate model set, we have
 \begin{align*}
  \sqrt {{T_2}} \left( {\hat \Delta _1(\hat{\bm w})- {\Delta _1}} \right)\mathop  \to \limits^d \mathcal{N}\left(0,\phi^2\bm{\mathcal{H}}\bm\Phi \bm\Xi \bm\Phi ^\top\bm{\mathcal{H}}^\top +\Sigma_v\right).
\end{align*}
\end{thm}

\begin{remark}\label{rk4}
Theorem \ref{th.4} confirms that {\sf RRMA} assigns all weights to the quasi-correct models asymptotically.
Moreover, it is easy to see that $ \hat w_{\rm{cor}}=\sum\limits_{s\in \mathcal{I}_{\rm{cor}}} {\hat w_s}\mathop  \to \limits^p  \sum\limits_{s\in \mathcal{I}_{\rm{cor}}} {w_s^*}=1$ by Theorems \ref{th.2} and \ref{th.4}, which implies $w_s^*=0$ for $s\notin \mathcal{I}_{\rm{cor}}$ due to $\bm w^*\in\mathcal{W}$. Hence,
we can deduce that $\Big(\sum_{s\notin \mathcal{I}_{\rm{cor}}} {\hat w_s}\Big)^2= \Big(\sum_{s\notin \mathcal{I}_{\rm{cor}}} {\left(\hat w_s-w_s^*\right)}\Big)^2
\leq \left(S-\left|\mathcal{I}_{\rm{cor}}\right| \right)\sum_{s\notin \mathcal{I}_{\rm{cor}}} {\left(\hat w_s-w_s^*\right)^2}
\leq S \left\| {\bm{\hat w} - \bm w^*} \right\|_2^2
=O_p\big( S^2\bar p^4\log T_1/T_1  \big)$, which means that the weight on any under-fitted candidate model is asymptotically zero.
If only one model among the candidate set is correctly specified, {\sf RRMA} will asymptotically pick that model.
When one or more quasi-correct models are included in the candidate model set, Theorem \ref{th.5} ensures that the proposed {\sf RRMA} converges to the true ATE at the rate of $ S^2\bar p^4\log T_1/T_1 +T_2^{-1}$.
Theorem \ref{th.6*} indicates the bias will disappear and $\hat\Delta_1(\hat {\bm w})$ reduces to an asymptotically unbiased estimator if there is at least one quasi-correct model in the candidate model set.

\end{remark}

\section{Implementation details}\label{sect.4}
\subsection{Prepare the non-nested candidate models}\label{sect.4A}

There are $2^{N-1}-1$ non-empty candidate models
for $N-1$ control units, involving an enormous number of possibilities for even just a moderate $N$.
Conducting an exhaustive model averaging over all possible models is computationally infeasible. Therefore, it becomes necessary to carefully construct a manageable and representative set of candidate models. Most existing works \citep{hr12,LS15,FYT23,JLLC24} considered using nested candidate models with a total of $N-1$ models. However, the $(s+1)$-th candidate model obviously includes repetitive information from the previous $s$ nested models, potentially inflating the variability and instability for the estimated weights. It is also unclear how to derive an optimal model sequence for such nested models.

We introduce a different yet attractive scheme to build a series of non-nested candidate models that can avoid overlapped information from different candidate models. Specifically, to safeguard against outliers and extreme observations, we first calculate the $p$-value between each covariate and the response by using the Kendall's rank correlation test \citep{LPZZ12}. Then we rank the $N-1$ covariates based on their $p$-values from low to high, and discard covariates whose $p$-values are above 5\%. Specifically, let $p^\natural_j$ be the $p$-value corresponding to the $j$-th covariate, $j=1,\cdots, N-1$. We write the index set of significant covariates as
\begin{align} \label{eq7B}
\mathcal{A}^\natural=\{j_1,\cdots,j_d: p^\natural_{j_1}\leq \cdots \leq p^\natural_{j_d}  \leq 5\%\} .
\end{align}

Next we partition $\mathcal{A}^\natural$ into $S$ mutually exclusive groups
\begin{align} \label{eq7c}
\overbrace {j_1,\cdots,j_m}^{\mathcal{A}_1},\overbrace {j_{m+1},\cdots,j_{2m}}^{\mathcal{A}_2},\cdots,\overbrace {j_{(S-2) m+1} ,\cdots,j_{(S-1) m}}^{\mathcal{A}_{S-1}},\overbrace {j_{(S-1)m+1},\cdots,j_{d}}^{\mathcal{A}_S}.
\end{align}
Each group is used once for constructing a candidate model, yielding the $S$ candidate models $\mathcal{M}_1,\cdots,\mathcal{M}_{S}$. Thus, the number of covariates in the model $\mathcal{M}_s$ is $|\mathcal{A}_s|=m$ for $s=1,\cdots,S-1$ and the number of covariates in the last model $\mathcal{M}_S$ is $|\mathcal{A}_S|=d-(S-1) \times m $. Set $S=\lfloor d/m\rfloor $, where $\lfloor a\rfloor$ is the largest integer less than or equal to $a$.
Note that only covariates without a significant marginal rank correlation are excluded at the sorting stage.

\begin{algorithm}[H]
    \caption{: A candidate model determination algorithm.
    }
    \label{alg:1}
    \begin{algorithmic}[1]
     \footnotesize

  \STATE Input: $\big\{ Y_{1t},\bm X_{t}\big\}_{t=1}^{T_1} $, $M_{\min}$ and $M_{\max}$.
 \STATE Obtain the initial significant index set $\mathcal{A}^\natural$ by (\ref{eq7B}).
  \FOR{$t=\frac{T_1}{2}+1,\frac{T_1}{2}+2,\cdots, T_1$,}
    \FOR{$m=M_{\min},M_{\min}+1,\cdots, M_{\max}$,}
   \STATE Set $S(m)=\lfloor d/m\rfloor$.
       \FOR{$s=1,2,\cdots, S(m)$,}
            \STATE  Set $\mathcal{A}_1,\cdots, \mathcal{A}_{S(m)}$ by replacing $S$ in (\ref{eq7c}) with $S(m)$. Calculate $\hat \mu_{1[t]}^{(s)}=\hat\alpha_{[t]}^{(s)}+\bm X_t^{(s)\top}\hat{\bm\beta}_{[t]}^{(s)} $, where $\hat{\bm\beta}^{(s)}_{[t]} $ and $\hat\alpha_{[t]}^{(s)}$ are obtained by \eqref{eq0p} and $\bm X_{t}^{(s)}=\big(X_{jt}: j\in \mathcal{A}_s\big)^\top$.
        \ENDFOR
        \STATE Calculate the AAPE: $\mathcal{Q}(m)= \frac{2}{T_1}\sum_{t=\frac{T_1}{2}+1}^{T_1}\left| Y_{1t}-\left(\hat{ \bm\mu}_{1[t]}^m\right)^{\top}\hat{\bm w}(m)\right|$, where $\hat{\bm w}(m)$ is obtained in Section \ref{sect.2.3} and
$\hat{ \bm\mu}_{1[t]}^m=\left(\hat \mu_{1[t]}^{(1)},\cdots,\hat \mu_{1[t]}^{(S(m))} \right)^\top$.
    \ENDFOR
  \ENDFOR
   \STATE  The optimal $m$ is determined by $\hat {m}= \mathop {\arg \min }\limits_{m\in \{M_{\min},\cdots, M_{\max}\}} \mathcal{Q}(m)$. Set $\hat S(\hat m)=\lfloor d/\hat {m}\rfloor$.
   \STATE Output: the final index sets $\mathcal{A}_1,\cdots, \mathcal{A}_S$ are obtained by replacing $m$ and $S$ in (\ref{eq7c}) with $\hat m$ and $\hat S(\hat m)$ respectively.
    \end{algorithmic}
\end{algorithm}

 Algorithm \ref{alg:1} groups the covariates with similar importance together to form candidate models. The individual model prediction power is expected to increase with the segment length or the sub-model size $m$. Similar to the scheme of choosing model weights in Section \ref{sect.2.3}, we use the GWCV to select the optimal $m$ by minimizing the average absolute prediction error (AAPE). We give the details on how to determine the index sets in Algorithm \ref{alg:1}.

\begin{remark}\label{rk2}
Algorithm \ref{alg:1} provides us with a new scheme to build a series of non-nested sub-models, which effectively balances estimation accuracy against computational complexity through systematic trade-off optimization. In addition, we notice that the number of candidate models $S $ and the model weight vector $\bm w$ depend on $m$. Thus, we rewrite $S$ and $\bm w$ as $S(m)$ and $\bm w(m)$ respectively in Algorithm \ref{alg:1}. In other places, we still use $S$ and $\bm w$ for the simplicity of notation. We take $M_{\min}=3$ and $M_{\max}=10$ in simulations, which lead to quite satisfactory estimation performance for the ATE, see Section \ref{sect.5}.
\end{remark}

\subsection{A confidence interval for the ATE}\label{sect.4B}

As shown in Theorem \ref{th.3} of Section \ref{sect.3.2}, the {\sf RRMA} estimator $\hat \Delta _1(\hat {\bm w})$ asymptotically follows a normal distribution, but the asymptotic variance $\sqrt{T_2}\hat \Delta _1(\hat {\bm w})$ involves the unknown matrices $\bm\Phi \bm\Xi \bm\Phi^\top$ (asymptotic covariance matrix of $\sqrt{T_1}\hat{\bm \theta}$) and $\Sigma_v$. It is challenging to estimate those covariance parameters accurately.
As an alternative, bootstrap methodologies are very useful in practice. However, the traditional i.i.d. bootstrap might be invalid in presence of time dependence. Thus, we utilize the well-known moving block bootstrap (MBB) proposed by \cite{Kunsch1989} to estimate the asymptotic covariance matrix of $\sqrt{T_1}\hat{\bm \theta}$.

We partition the original samples $\big\{ Y_{1t},\bm X_{t}\big\}_{t=1}^{T_1} $ into overlapping blocks of length $l$, denoted as $\mathbb{B}_1,\cdots,\mathbb{B}_{L}$, where $\mathbb{B}_t=\left\{ \left(Y_{1t},\bm X_{t}\right),\left(Y_{1t+1},\bm X_{t+1}\right),\cdots,\left(Y_{1t+l-1},\bm X_{t+l-1}\right)\right\}$ for $t=1,\cdots,L$ and $L = T_1-l+1$ is the number of blocks.
Set $l=\lfloor T_1^{1/3} \rfloor$ and $k=\lceil T_1/l\rceil$, where $\lfloor a\rfloor$ ($\lceil a \rceil$) is the largest (smallest) integer less (greater) than or equal to $a$.
Then, we select $k$ blocks at random with replacement from the
collection $\mathbb{B}_1,\cdots,\mathbb{B}_{L}$, denoted as $\mathbb{B}_{1}^\dag,\cdots,\mathbb{B}_{k}^\dag$.
Join the $k$ blocks in the order $\mathbb{B}_{1}^\dag,\cdots,\mathbb{B}_{k}^\dag$ together to obtain a new
set of the bootstrap samples of length $T_1$ denoted by $\big\{ Y_{1t}^\dag,\bm X_{t}^\dag\big\}_{t=1}^{T_1} $.
Then, we use the same calculation procedure of Section \ref{sect.2.2} to obtain the bootstrap estimator $\hat{\bm \theta}^\dag$ based on $\big\{ Y_{1t}^\dag,\bm X_{t}^\dag\big\}_{t=1}^{T_1} $.
We repeat the above process for a large number of times (i.e., $B$ times), and obtain a series of bootstrap estimators $\big\{ \hat{\bm \theta}_b^\dag  \big\}_{b=1}^{B}$. Then we can estimate the asymptotic covariance matrix of $\sqrt{T_1}\hat{\bm \theta}$ by
$
\widehat{\rm Cov}\left(\sqrt{T_1}\hat{\bm \theta}\right)=\frac{T_1}{B-1}\sum\limits_{b = 1}^{B}{\left(\hat{\bm\theta}_b^\dag -\bar{\hat{\bm\theta}}_b^\dag\right)\left(\hat{\bm\theta}_b^\dag -\bar{\hat{\bm\theta}}_b^\dag\right)^\top},
$
where $\bar{\hat{\bm\theta}}_b^\dag=\frac{1}{B}\sum\limits_{b = 1}^{B}{\hat{\bm\theta}_b^\dag}$. We set $B=1000$ in numerical studies.

Following \cite{L20}, we estimate the variance $\Sigma_v$ by
$\widehat{\Sigma}_v=\frac{1}{T_2}\sum\limits_{t = T_1+1}^{T}{{\hat v_t}^2}$,
where $\hat v_t=\hat \Delta _{1t}(\hat {\bm w})-\hat \Delta _1(\hat {\bm w})$ with $\hat \Delta _{1t}(\hat {\bm w})$ and $\hat \Delta _1(\hat {\bm w})$ specified in \eqref{eq7}.
Then, combined with the result of Theorem \ref{th.3}, the asymptotic variance of $\hat \Delta _1(\hat {\bm w})$ can be estimated by
\begin{align*}
\widehat{\rm Var}\left(\hat \Delta _1(\hat {\bm w})\right)=\frac{1}{T_2}\hat\phi^2\hat{\bm{\mathcal{H}}}\widehat{\rm Cov}\left(\sqrt{T_1}\hat{\bm \theta}\right)\hat{\bm{\mathcal{H}}}^\top+\frac{1}{T_2} \widehat{\Sigma}_v,
\end{align*}
where $\hat\phi =\sqrt {{T_2}/{T_1}}$, $\hat{\bm{\mathcal{H}}}=\left( \hat{\bm{\mathcal{H}}}^{(1)},\cdots, \hat{\bm{\mathcal{H}}}^{(S)}\right)$,
$\hat{\bm{\mathcal{H}}}^{(s)}= \left(\hat w_s, \hat w_sT_2^{ - 1}\sum\limits_{t = {T_1} + 1}^T \bm  X_t^{(s)\top}\right)$.
Thus, an asymptotic $1-\alpha$ confidence interval for $\Delta_1$ is
\begin{align}\label{eqCI}
{\rm CI}_{1-\alpha}=\left[\hat \Delta _1(\hat {\bm w})-z_{\alpha/2}\sqrt{\widehat{\rm Var}\left(\hat \Delta _1(\hat {\bm w})\right)},\hat \Delta _1(\hat {\bm w})+ z_{\alpha/2}\sqrt{\widehat{\rm Var}\left(\hat \Delta _1(\hat {\bm w})\right)} \right],
\end{align}
where $z_{\alpha/2}$  is the upper $\alpha/2$ quantile of the standard normal distribution.

\section{Simulation studies}\label{sect.5}
In this section, we will investigate the finite sample performance of the proposed model averaging estimator. We conduct numerical experiments with R running on a desktop computer with a CPU specification (11th Gen Intel(R) Core(TM) i7-11700, 2.50GHz, 16.00 GB of RAM). To make a full comparison, we consider the following competitive approaches.

(a) {\sf PHCW}: \cite{LB17} suggested to estimate $\alpha$ and $\bm\beta$ in the full model \eqref{eqA} by the penalized least squares method, that is $\big(\hat\alpha_{{\rm PHCW}},\hat{\bm\beta}_{{\rm PHCW}}^\top\big)=\mbox{arg} \mathop {\min }\limits_{\alpha \in \mathbb{R},~\bm\beta \in \mathbb{R}^{N-1}} \sum_{t=1}^{T_1}\big(Y_{1t}-\alpha-\bm X_t^\top\bm\beta\big)^2+\sum_{j=1}^{N-1}p_{\lambda}\big(|\beta_j|\big)$, where $p_{\lambda}(\cdot)$ is taken as the SCAD penalty and the tuning parameter $\lambda$ is selected by 10-fold cross validation. Then, the counterfactual outcome $Y_{1t}(0)$ is estimated by $\hat \mu_{1t}= \hat\alpha_{{\rm PHCW}}+\bm X_t^\top\hat{\bm\beta}_{{\rm PHCW}}$ and the final ATE estimator is obtained by (\ref{eqA0}).

(b) {\sf MSCM}: \cite{L20} proposed to estimate $\alpha$ and $\bm\beta$ in the full model \eqref{eqA} by the modified synthetic control method, that is $\big(\hat\alpha_{{\rm MSCM}},\hat{\bm\beta}_{{\rm MSCM}}^\top\big)=\mbox{arg} \mathop {\min }\limits_{\alpha \in \mathbb{R},~\bm\beta \in \Lambda} \sum_{t=1}^{T_1}\big(Y_{1t}-\alpha-\bm X_t^\top\bm\beta\big)^2$, where $\Lambda=\left\{ \bm \beta\in \mathbb{R}^{N-1}: \beta_j\geq 0~\mbox{for}~ j=1,\cdots,N-1\right\}$. Then, the counterfactual outcome $Y_{1t}(0)$ is estimated by $\hat\alpha_{{\rm MSCM}}+\bm X_t^\top\hat{\bm\beta}_{{\rm MSCM}}$ and the final ATE estimator is obtained by (\ref{eqA0}).

(c) {\sf LSMA}: \cite{GLW15} introduced a nested model averaging strategy for estimating the ATE, where the parameters in each candidate model are estimated by the least squares method and the model weights are determined by the least
squares cross-validation criterion.

(d) {\sf EW}, {\sf BIC}, {\sf RRMA}: It is known that the performance of model averaging heavily relies on the weight selection criterion, and thus we consider different approaches for selection of model weights. {\sf EW} simply assigns equal weights to all candidate models. {\sf SBIC} assigns the model weights by the smoothed BIC criterion, implying that the model weight of the $s$ candidate model is determined by ${{\hat w}_s}
=\exp \left( { - 0.5{\rm {BIC}}_{(s)}} \right)/
\sum\limits_{j = 1}^S {\exp \left( { - 0.5{\rm {BIC}}_{(j)}} \right)}$ with ${\rm {BIC}}_{(s)}={\rm log}\left( \frac{1}{{{T_1}\left( {{T_1} - 1} \right)}}\sum_{t =1}^{T_1}\sum_{t' \neq t} {b_{tt'}^{(s)}\left| Y_{1tt'}-\bm X_{tt'}^{(s)\top}\bm\beta^{(s)}\right|} \right) +p_s\log T_1 /T_1 $. {\sf RRMA} is the proposed optimal model averaging.

In simulation studies, we calculated the bias, the mean absolute error (MAE), the mean squared error (MSE), the coverage probabiity (CP) and the average length (AL) of confidence intervals to assess estimation accuracy for the ATE. The number of Monte Carlo simulation replications is $500$. Note that {\sf PHCW} and {\sf LSMA} do not provide the confidence interval for the ATE, and hence we do not report CP and AL for them.

\textbf{Example 1} Consider the stationary three-factor model. The samples of the potential outcomes are generated by the model (\ref{eqfoc*}) where $\mathbb{K}=3$ and $\bm B=\left(\bm b_1,\cdots,\bm b_N\right)^\top$ is the $N\times 3$ loading matrix with the $i$-th row $\bm b_i=\left(b_{1i},b_{2i},b_{3i}\right)^\top$ and $\bm f_t=\left(f_{1t},f_{2t},f_{3t}\right)^\top$ is the three-dimensional common factor vector. We set $a_1=\cdots=a_N=2$, $f_{1t}=0.8f_{1t-1}+\upsilon_{1t}$, $f_{2t}=-0.6f_{1t-1}+\upsilon_{2t}+0.8\upsilon_{2t-1}$, $f_{3t}=\upsilon_{3t}+0.9\upsilon_{3t-1}+0.4\upsilon_{3t-2}$ and $\upsilon_{jt}\mathop\sim\limits^{iid}\mathcal{N}\left(0,0.5^2\right)$ for $j=1,2,3$. The factor loading $ b_{ji}$ is taken by
\[b_{ji}=\left\{ {\begin{array}{*{20}{c}}
   {2,~i\in \{1,11,\cdots,21\}}  \\
   {b_{ji}^\dag,~i\in \{2,3,\cdots,10\}}  \\
   b_{ji}^\ddag,~i\in \{22,23,\cdots,N\}  \\
\end{array}} \right. j=1,2,3,\]
where $b_{ji}^\dag\mathop\sim\limits^{iid}\mathcal{N}\left(0,1\right)$ and
$b_{ji}^\ddag\mathop\sim\limits^{iid}\mathcal{N}\left(0,0.5^2\right)$.
We note that the first unit (the response) and the $j$-th control unit are positively correlated with shared common factors for $j\in \{11,\cdots,21\}$ due to large positive loadings, and the remaining control units may be weakly correlated.  To incorporate the treatment effect on the first unit, we set $Y_{1t}(1)=Y_{1t}(0)+\Delta_{1t}$ for the post-treatment time points $t=T_1+1,\cdots,T$, where the individual treatment effect $\Delta_{1t}$ is generated by $\Delta_{1t}=\alpha_0\left(\exp(z_t)/(1+\exp(z_t))+1\right)$ with $z_t=0.5z_{t-1}+\iota_t$ and $\iota_t\mathop\sim\limits^{iid}\mathcal{N}\left(0,0.5^2\right)$. Then, the true ATE is constructed by $\Delta_{1}=T_2^{-1}\sum\limits_{t=T_1+1}^{T}\Delta_{1t}$. Clearly, the treatment effect is positive when $\alpha_0>0$, and there is no effect when $\alpha_0=0$.

To study the robustness of the estimators, three distributions of the random error $u_{it}$ are considered:
\emph{Case a}: $u_{it}\mathop\sim\limits^{iid}\mathcal{N}\left(0,1\right)$ for $1\leq t\leq T$; \emph{Case b}: $u_{it}\mathop\sim\limits^{iid}t_3$ for $1\leq t\leq T_1$ and $u_{it}\mathop\sim\limits^{iid}\mathcal{N}\left(0,1\right)$ for $T_1+1\leq t\leq T$;
\emph{Case c}: $u_{it}\mathop\sim\limits^{iid}\mathcal{CN}$ for $1\leq t\leq T_1$ and $u_{it}\mathop\sim\limits^{iid}\mathcal{N}\left(0,1\right)$ for $T_1+1\leq t\leq T$;
where $t_3$ and $\mathcal{CN}$ stand for the $t$-distribution with three degrees of freedom and the mixture of $\mathcal{N}\left(0,1\right)$ and $\mathcal{N}\left(0,25\right)$ with weights 90$\%$ and 10$\%$, respectively.

Table \ref{table1AA} report the estimation results of different approaches.
First, although all candidate models may be mis-specified in practice, {\sf RRMA} yields relative small biases when compared with other existing approaches.
Second, it is easy to see that {\sf RRMA} outperforms {\sf EW} and {\sf SBIC} under all situations, improving estimation accuracy by selecting better model weights.
Third, the traditional single model-based methods ({\sf PHCW} and {\sf MSCM}) are inferior to the proposed model averaging procedure {\sf RRMA} in all cases. One possible explanation is that {\sf PHCW} and {\sf MSCM} may suffer from serious model mis-specification. Fourth, for the normal error (\emph{Case a}), {\sf RRMA} performs comparably to
{\sf LSMA}, {\sf PHCW} and {\sf MSCM}, but it performs best among all approaches for non-normal errors (\emph{Case b} and \emph{Case c}). Finally, the confidence intervals for {\sf RRMA} and {\sf MSCM} work reasonably well since the coverage probabilities of {\sf RRMA} and {\sf MSCM} are close to the nominal value, but the remaining methods ({\sf EW}, {\sf SBIC}) significantly underestimate the coverage probabilities for some cases (e.g., $N=45,\alpha_0=0$). Furthermore, the average length of the confidence interval of {\sf RRMA} is shorter than that of {\sf MSCM}, especially for \emph{Case b} and \emph{Case c}.

Figure S2 of the Supplementary Material displays the probability density function of {\sf RRMA} based on 500 repetitions for $T_1=60,T_2=30,N=25$.
We find that the blue dashed line is close the black solid one and the $p$-values of Shapiro-Wilk normality test are always greater than 0.05, which shows that {\sf RRMA} follows closely a normal distribution whether a quasi-correct model is included in the candidate model set or not. Hence, the normal approximation confidence interval given in \eqref{eqCI} works in practice.

\tabcolsep=4.5pt
\begin{table}
\caption{\small The bias, mean absolute error (MAE) and mean squared error (MSE $\times 100$) of various methods among 500 replications for Example 1. }
\label{table1AA}\small
\begin{tabular}{cccccccccccccc} \noalign{\smallskip}\hline
\multirow{2}{*}{$error$}
&\multirow{2}{*}{method}
&\multicolumn{5}{c}{$\alpha_0=0$}
&&\multicolumn{5}{c}{$\alpha_0=1$}
\\
\cline{3-7}
\cline{9-13}
&&\multirow{1}{*}{{\rm Bias}}&\multirow{1}{*}{{\rm MAE}}&\multirow{1}{*}{{\rm MSE}}&\multirow{1}{*}{{\rm CP}}&\multirow{1}{*}{{\rm AL}}
&&\multirow{1}{*}{{\rm Bias}}&\multirow{1}{*}{{\rm MAE}}&\multirow{1}{*}{{\rm MSE}}&\multirow{1}{*}{{\rm CP}}&\multirow{1}{*}{{\rm AL}}
\\
\hline
&&\multicolumn{11}{c}{$T_1=60,T_2=30,N=25$}\\
$Case~ a$&{\sf RRMA}&-0.010 &0.227&  7.568& 0.960& 1.192&&-0.003& 0.213&  7.165& 0.954& 1.246 \\
&{\sf EW}& -0.011& 0.258& 10.15& 0.942& 1.246&&-0.005& 0.253& 10.661& 0.924& 1.296\\
&{\sf SBIC}&-0.012& 0.247&  9.224& 0.946& 1.208&&-0.005& 0.239&  9.478& 0.928& 1.258\\
&{\sf PHCW}&-0.021& 0.237&  9.244& --&--&&0.009& 0.234&  9.042&--&--\\
&{\sf MSCM}&-0.015& 0.230&  8.301& 0.948& 1.143&&0.006& 0.215&  7.231& 0.956& 1.152\\
&{\sf LSMA}&-0.020& 0.246&  9.254&--&--&&0.014& 0.234&  8.666 &--&--\\

$Case~ b$&{\sf RRMA}&0.002& 0.258& 10.43& 0.950& 1.386&&-0.001& 0.247&  9.528& 0.952& 1.315\\
&{\sf EW}&0.001& 0.301& 14.281& 0.914& 1.335&&0.004& 0.290& 12.89& 0.932& 1.359\\
&{\sf SBIC}&0.001& 0.292& 13.44& 0.918& 1.312&&0.003& 0.281& 12.15& 0.926& 1.337\\
&{\sf PHCW}&-0.001& 0.339& 17.78&--&--&&0.009& 0.316& 15.90&--&--\\
&{\sf MSCM}&0.014& 0.289& 13.08& 0.948& 1.461&&0.014& 0.303& 15.01& 0.948& 1.499\\
&{\sf LSMA}&0.001& 0.324& 16.67&--&--&&0.018& 0.334& 18.17&--&--\\

$Case~ c$&{\sf RRMA}&0.008& 0.227&  8.100& 0.958& 1.266&&-0.003& 0.230&  8.316& 0.952& 1.494\\
&{\sf EW}&0.010& 0.270& 11.64& 0.940& 1.320&&0.003& 0.255& 10.40& 0.938& 1.471 \\
&{\sf SBIC}&0.010& 0.258& 10.59& 0.942& 1.290&&0.001& 0.244&  9.529& 0.942& 1.448\\
&{\sf PHCW}&0.015& 0.295& 14.25&--&--&&0.013& 0.307& 15.44&--&--\\
&{\sf MSCM}&-0.008& 0.258& 10.66& 0.924& 1.319&&0.005& 0.277& 12.27& 0.942& 1.344\\
&{\sf LSMA}&0.008& 0.309& 15.70&--&--&&0.018& 0.311& 15.23&--&--\\

&&\multicolumn{11}{c}{$T_1=120,T_2=60,N=45$}\\
$Case~ a$&{\sf RRMA}&0.006& 0.151& 3.564& 0.952& 0.722&& -0.010& 0.155& 3.774& 0.952 &0.728\\
&{\sf EW}& 0.016& 0.249& 9.284& 0.894& 0.953&& 0.005 &0.238& 8.778 &0.878 &0.949 \\
&{\sf SBIC}&0.015& 0.230& 8.009& 0.906& 0.900&&0.003 &0.222 &7.600 &0.886 &0.897\\
&{\sf PHCW}&0.010& 0.159& 4.052& --&--&& -0.011 &0.157 &3.943   & --&--\\
&{\sf MSCM}&0.008& 0.150& 3.574& 0.956& 0.791&& -0.013& 0.152& 3.751 &0.964& 0.789\\
&{\sf LSMA}&0.014& 0.156& 3.855& --&--&&-0.005 &0.154 &3.769   &--&--\\

$Case~ b$&{\sf RRMA}& 0.004& 0.177&  5.201& 0.946& 0.911&&-0.002& 0.177&  4.908 &0.946 &0.905\\
&{\sf EW}& 0.001& 0.254& 10.67& 0.884& 1.023&&-0.007 &0.264 &10.557 &0.876& 1.013\\
&{\sf SBIC}&0.001& 0.243&  9.773& 0.888& 0.989&&-0.007& 0.253 & 9.709 &0.878 &0.981\\
&{\sf PHCW}&0.004& 0.217&  7.832& --&--&&-0.002 &0.222 & 8.073  &  -- & --\\
&{\sf MSCM}&-0.011& 0.212&  7.003& 0.952& 1.013&&-0.005& 0.195 & 6.268 &0.950 &1.006\\
&{\sf LSMA}&-0.006& 0.218&  7.857& --&--&& 0.000 &0.208  &7.228  &  -- &--\\

$Case~ c$&{\sf RRMA}& -0.002& 0.145& 3.479& 0.948& 0.761&&0.007 &0.150& 3.574& 0.962 &0.762\\
&{\sf EW}& 0.011& 0.227& 8.283& 0.906& 0.950&&0.020 &0.238 &8.869 &0.884 &0.949\\
&{\sf SBIC}&0.009& 0.213& 7.267& 0.910& 0.906&&0.019& 0.223& 7.782 &0.892& 0.906\\
&{\sf PHCW}&-0.002& 0.202& 6.867& --&--&&0.012&0.209 &7.051  & -- &--\\
&{\sf MSCM}&0.001& 0.175& 4.935& 0.940& 0.929&&0.007& 0.184& 5.539& 0.904 &0.923\\
&{\sf LSMA}&0.003& 0.192& 6.089& --&--&&0.023 &0.202 &6.695  &-- &--\\
\hline
\end{tabular}
\\
Note: The symbol ``--" indicates that the result of a method is not available.
\end{table}

\tabcolsep=3.5pt
\begin{table}
\caption{\small The mean and standard deviation (in the parenthese) of the first component of {\sf RRMA} ($\hat w_1$) for Example 2, when the first candidate model is the quasi-correct model. }
\label{tableB}\small
\begin{tabular}{cccccccccc} \noalign{\smallskip}\hline
\multirow{2}{*}{$T_1$}
&\multicolumn{2}{c}{$Case~ a$}
&&\multicolumn{2}{c}{$Case~ b$}
&&\multicolumn{2}{c}{$Case~ c$}
\\
\cline{2-3}
\cline{5-6}
\cline{8-9}
&\multirow{1}{*}{$m=3$}&\multirow{1}{*}{$m=6$}
&&\multirow{1}{*}{$m=3$}&\multirow{1}{*}{$m=6$}
&&\multirow{1}{*}{$m=3$}&\multirow{1}{*}{$m=6$}
\\
\hline
&\multicolumn{8}{c}{$\alpha_0=0$}\\
60&0.917 (0.086) &0.906 (0.092)&&0.852 (0.148) &0.816 (0.166)&&0.920 (0.095) &0.888 (0.105) \\
80&0.946 (0.061)& 0.942 (0.062)&&0.879 (0.114)& 0.860 (0.112)&&0.930 (0.074)& 0.907 (0.085)\\
120& 0.954 (0.049)& 0.951 (0.050)&& 0.908 (0.091)& 0.901 (0.091)&&0.953 (0.055)& 0.935 (0.064)
\\
&\multicolumn{8}{c}{$\alpha_0=1$}\\
60& 0.911 (0.092)& 0.901 (0.091)&&0.861 (0.130)& 0.825 (0.135) &&0.904 (0.102)& 0.884 (0.115)\\
80& 0.933 (0.082)& 0.925 (0.087)&& 0.886 (0.110)& 0.865 (0.122)&&0.935 (0.078)& 0.905 (0.094)\\
120&0.962 (0.047)& 0.956 (0.049)&&0.910 (0.087)& 0.901 (0.090)&&0.951 (0.059)& 0.934 (0.064)\\
\hline
\end{tabular}
\end{table}

\tabcolsep=7pt
\begin{table}
\caption{\small The bias, mean absolute error (MAE), mean squared error (MSE), coverage probabiity (CP) and average length (AL) of {\sf RRMA}  among 500 replications for Example 2, when the first candidate model is the quasi-correct model.}
\label{table1C}\small
\begin{tabular}{cccccccccccc} \noalign{\smallskip}\hline
\multirow{2}{*}{$T_1,T_2$}
&\multicolumn{5}{c}{$m=3$}
&&\multicolumn{5}{c}{$m=6$}
\\
\cline{2-6}
\cline{8-12}
&\multirow{1}{*}{{\rm Bias}}&\multirow{1}{*}{{\rm MAE}}&\multirow{1}{*}{{\rm MSE}}&\multirow{1}{*}{{\rm CP}}&\multirow{1}{*}{{\rm AL}}
&&\multirow{1}{*}{{\rm Bias}}&\multirow{1}{*}{{\rm MAE}}&\multirow{1}{*}{{\rm MSE}}&\multirow{1}{*}{{\rm CP}}&\multirow{1}{*}{{\rm AL}}
\\
\hline
&\multicolumn{11}{c}{$Case~ a$, $\alpha_0=0$}\\
60,30& -0.009 & 0.259 & 0.104&  0.936&  1.216&& -0.002 & 0.265 & 0.112&  0.948&  1.375\\
80,40& 0.019& 0.223& 0.080& 0.946& 1.014&& 0.014& 0.225& 0.079& 0.946& 1.072\\
120,60&0.006 &0.167 &0.042 &0.964 &0.814 &&0.004& 0.169 &0.043 &0.966 &0.842\\

&\multicolumn{11}{c}{$Case~ b$, $\alpha_0=0$}\\
60,30&-0.032&  0.278&  0.120&  0.948&  1.409&& -0.034&  0.296&  0.139 & 0.964&  1.721\\
80,40&0.009& 0.242& 0.092& 0.954& 1.190&& 0.007& 0.260& 0.101& 0.956& 1.285\\
120,60&-0.021  &0.188 & 0.059 & 0.950 & 0.943 &&-0.030 & 0.192&  0.062 & 0.956  &0.986
\\

&\multicolumn{11}{c}{$Case~ c$, $\alpha_0=0$}\\
60,30&-0.019&  0.279&  0.121 & 0.936 & 1.296 &&-0.019 & 0.288 & 0.135&  0.938&  1.454\\
80,40& -0.012 & 0.235&  0.084&  0.946&  1.089&& -0.018&  0.247&  0.089&  0.956&  1.167\\
120,60& -0.010  &0.176 & 0.049  &0.950 & 0.864 &&-0.011&  0.181 & 0.052 & 0.960 & 0.904\\
\hline
&\multicolumn{11}{c}{$Case~ a$, $\alpha_0=1$}\\
60,30&-0.019&  0.256&  0.110&  0.942&  1.222&& -0.029 & 0.265 & 0.108&  0.960 & 1.397\\
80,40&0.007& 0.209& 0.072& 0.948& 1.031&& 0.009& 0.225& 0.081& 0.944& 1.109\\
120,60&0.023 &0.166& 0.042& 0.944& 0.817&& 0.028& 0.160& 0.040& 0.970& 0.854\\

&\multicolumn{11}{c}{$Case~ b$, $\alpha_0=1$}\\
60,30& -0.017&  0.278&  0.122&  0.952&  1.415&& -0.020&  0.286&  0.132&  0.970&  1.613\\
80,40& 0.018& 0.255& 0.102& 0.946& 1.193&& 0.026& 0.260& 0.104& 0.950& 1.287\\
120,60&0.012& 0.193& 0.059& 0.958& 0.951&& 0.003 &0.194 &0.062& 0.962& 0.987\\

&\multicolumn{11}{c}{$Case~ c$, $\alpha_0=1$}\\
60,30& -0.016&  0.274&  0.119&  0.946&  1.287&& -0.018&  0.293&  0.143 & 0.946 & 1.438\\
80,40&-0.005 & 0.233&  0.087 & 0.956&  1.114&& -0.014 & 0.244&  0.097&  0.952&  1.190\\
120,60& -0.001&  0.171&  0.046&  0.952&  0.902&& -0.004&  0.173&  0.048 & 0.972&  0.962\\
\hline
\end{tabular}
\end{table}

\textbf{Example 2} The samples are still generated form the stationary three-factor model in Example 1. The factor loading is set as $ b_{j1}=2$, $b_{ji}=-2$ for $i=2,\cdots,4$ and $b_{ji}=0$ for $i=5,\cdots,N=15$. Now only the first four units are strongly correlated with the common factors with non-zero loadings, which indicates that only the first three covariates are important to predict the response. Figure S1 of the Supplementary Material displays boxplots of the Kendall rank correlation between the response and the covariates among 500 replications.

The goal of this example is to examine the situation when there is a quasi-correct model in the candidate model set. Although Algorithm \ref{alg:1} provides a reasonable set of candidate models, it does not automatically guarantee that a quasi-correct model always belong to the model set at each simulation. Thus, we directly partition the index set of $N-1$ predictors $\{1,\cdots,N-1\}$ into $S$ mutually exclusive groups $\mathcal{A}_1,\cdots,\mathcal{A}_S$, where $\mathcal{A}_s=\{(s-1)\times m+1,\cdots,s\times m\}$ for $s=1,\cdots,S-1$ and $\mathcal{A}_S=\{(S-1)\times m+1,\cdots,N-1\}$. Clearly, when $m\geq 3$, the first set includes all three important covariates, leading to a proper quasi-correct model.

Table \ref{tableB} reports means and standard deviations of the estimated weight for the quasi-correct model $\hat w_1$. As $T_1$ becomes large, $\hat w_1$ approaches one as expected.
Table \ref{table1C} presents the Bias, MAE, MSE, CP and AL for the proposed model averaging with $m=3$ and $6$. We observe that MAE and MSE decrease quickly and drop to zero as $T_1$ and $T_2$ increase. This verifies the consistency established in Theorem \ref{th.5}. Meanwhile, the proposed interval estimator also works quite well for all cases.

\section{An empirical analysis of terrorism's economic impacts in the Basque Country}\label{sect.6}

To empirically demonstrate the applicability of the proposed model averaging estimator,
we conduct a counterfactual analysis of the \emph{Basque Country}'s per capita GDP,  evaluating how its economy might have evolved in the absence of terrorist conflict.The analysis utilizes a region-level panel dataset originally studied by \citet{AG03}, which includes per capita GDP data for 17 Spanish regions spanning the years 1955 to 1997. This dataset is publicly available via the R package {\sf Synth}. The \emph{Basque Country} was among Spain most economically prosperous regions until the mid-1970s. However, the onset of widespread terrorist activities in the region during this period marked the beginning of a prolonged phase of economic decline. Political instability, particularly in the form of terrorism, is widely acknowledged as a major impediment to economic development. Figure S3 of the Supplementary Material presents the per capita GDP of 17 Spanish regions during two key intervals: 1955--1974 and 1990--1997. It is evident that the \emph{Basque Country} ranking declined from third to sixth place over a 15-year span characterized by persistent political unrest. Accordingly, we define the pre-treatment period as 1955--1974 and the post-treatment period as 1975--1997.
We take the \emph{Basque Country} as the treatment group, while the remaining 16 regions serve as control units, as they did not experience terrorist activity or political conflict during the study period. Specifically, the 16 control regions are:
\emph{Andalucia} ($X_{1t}$), \emph{Aragon} ($X_{2t}$), \emph{Asturias} ($X_{3t}$), \emph{Baleares}
 ($X_{4t}$), \emph{Canarias} ($X_{5t}$), \emph{Cantabria} ($X_{6t}$), \emph{Leon} ($X_{7t}$), \emph{Mancha} ($X_{8t}$),
\emph{Catalonia} ($X_{9t}$), \emph{Valenciana} ($X_{10t}$), \emph{Extremadura} ($X_{11t}$), \emph{Galicia} ($X_{12t}$), \emph{Madrid} ($X_{13t}$), \emph{Murcia} ($X_{14t}$), \emph{Navarra}($X_{15t}$) and \emph{Rioja} ($X_{16t}$), where $t=1955,\cdots,1997$.

\tabcolsep=2pt
\begin{table}\centering
\small
  \caption{{\small Confidence intervals and their lengths (Len.) (in parenthesis) of various approaches for the Basque Country dataset.}}
  \label{tab:CIs}
  \begin{tabular}{lcccc}
    \hline
          &90\%CI (Len.) &  95\%CI (Len.)&  99\%CI (Len.)\\
    \hline
{\sf RRMA}&$[-0.8027, -0.3925]\ (0.4102)$& $[-0.8420, -0.3532]\ (0.4888)$& $[-0.9188, -0.2764]\ (0.6424)$\\
{\sf EW}&$[-1.4614, -0.4893]\ (0.9721)$& $[-1.5545, -0.3962]\ (1.1583)$& $[-1.7365, -0.2142]\ (1.5223)$\\
{\sf SBIC}&$[-1.3193, -0.5102]\ (0.8091)$& $[-1.3968, -0.4327]\ (0.9641)$& $[-1.5483, -0.2812]\ (1.2671)$\\
{\sf MSCM}&$[-1.0638, -0.5413]\ (0.5225)$& $[-1.1083, -0.5089]\ (0.5994)$& $[-1.1806, -0.4325]\ (0.7481)$\\
    \hline
  \end{tabular}
\end{table}

The $p$-values between the \emph{Basque Country} and each control unit based on the Kendall rank correlation test are all less than $5\%$. Hence, no control units are removed in the initial step of Algorithm \ref{alg:1}. We take $M_{\min}=2$ and $M_{\max}=7$ in here.
Directly applying Algorithm \ref{alg:1} yields the optimal number of candidate models to be $S=8$, indicating that each candidate model includes two covariate, namely, $\mathcal{M}_1=\{X_{1t},X_{5t}\},\mathcal{M}_2=\{X_{2t},X_{3t}\},\mathcal{M}_3=\{X_{4t},X_{6t}\},
\mathcal{M}_4=\{X_{7t},X_{8t}\},\mathcal{M}_{5}=\{X_{9t},X_{10t}\},\mathcal{M}_{6}=\{X_{11t},X_{12t}\},
\mathcal{M}_{7}=\{X_{14t},X_{15t}\},
\mathcal{M}_{8}=\{X_{13t},X_{16t}\}$. Since {\sf PHCW} drops the non-negative and sum-to-one restrictions, some of the estimated weights are negative and greater than one. {\sf LSMA} estimates the ATE via a weighted average of 16 nested candidate models with covariates $\{X_{1t}\},\{X_{1t},X_{2t}\},\cdots,\{X_{1t},\cdots,X_{16t}\}$.

Figure S4 of the Supplementary Material displays the estimated coefficients for {\sf PHCW} and {\sf MSCM} and the estimated model weights for {\sf LSMA}, {\sf RRMA}, {\sf EW} and {\sf SBIC}. We notice that {\sf RRMA} assigns the positive weights to $\mathcal{M}_1$, $\mathcal{M}_4$, $\mathcal{M}_5$ and $\mathcal{M}_8$. Thus, the trends in per capita GDP growth of \emph{Basque Country} prior to terrorism is well modelled by these four models.
In particular, \emph{Catalonia} ($X_{9t}$) is assigned the greatest weight by {\sf MSCM} and {\sf LSMA}, and {\sf RRMA} assigns the relatively weight to the fifth model $\mathcal{M}_{5}$ which also includes $X_{9t}$.
Figure S5 (a) of the Supplementary Material clearly shows that the per capita GDP of \emph{Catalonia} is most similar to that of \emph{Basque Country} before the treatment.
Moreover, this result is also consistent with the findings of \cite{AG03}, who identified \emph{Catalonia} as exhibiting comparable socioeconomic attributes to \emph{Basque Country}, particularly in terms of industrial composition and political autonomy.
Figure S5 (b) provides various estimators for the ATE, and Table \ref{tab:CIs} reports the confidence intervals and the corresponding interval lengths. We notice that the upper bounds of almost all intervals are negative, suggesting that the estimated ATE values are significantly different from zero. Moreover, the interval lengths of {\sf RRMA} are the smallest among all competing methods.

\begin{figure}
\centering
\includegraphics[scale=0.2]{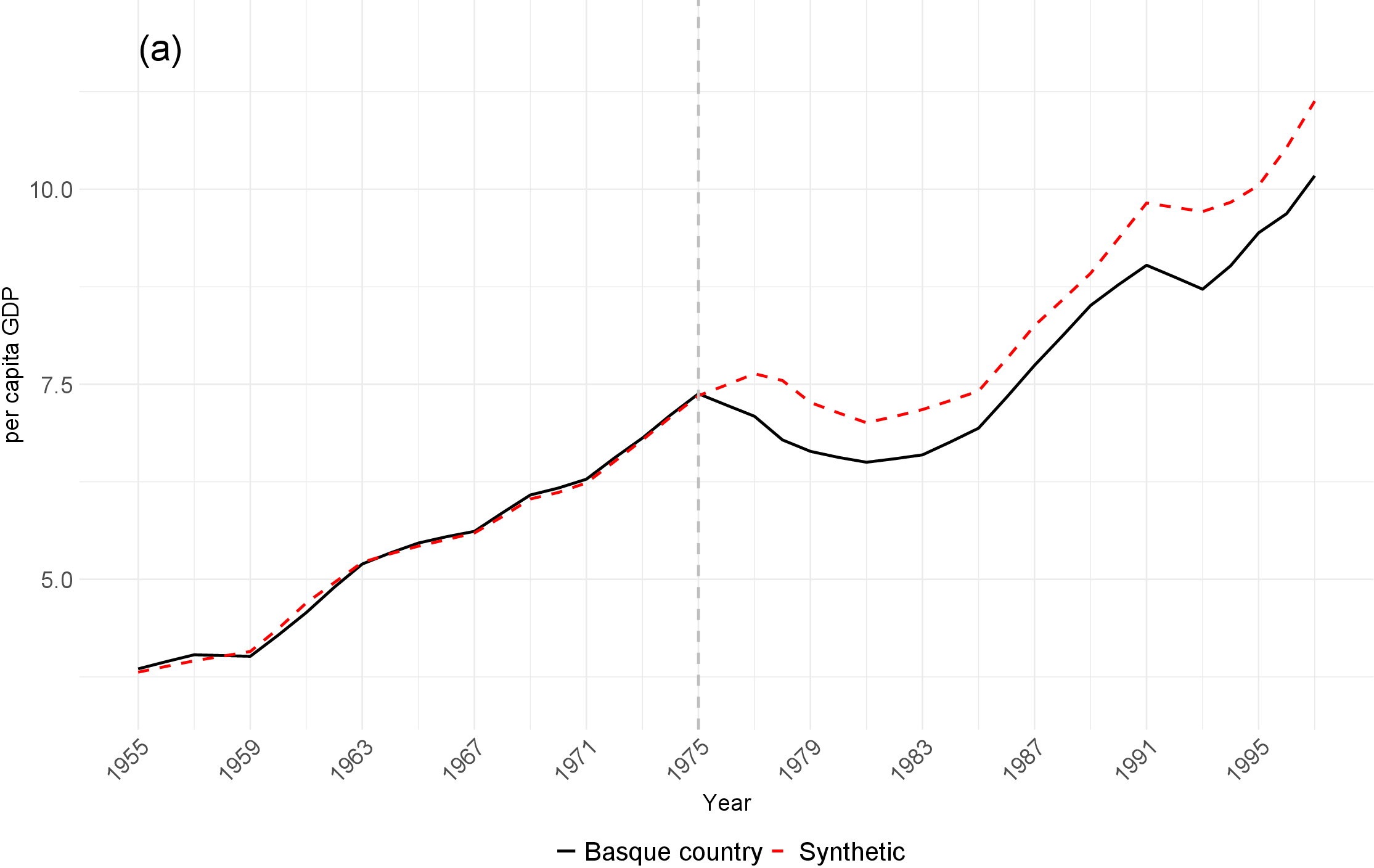}
\includegraphics[scale=0.2]{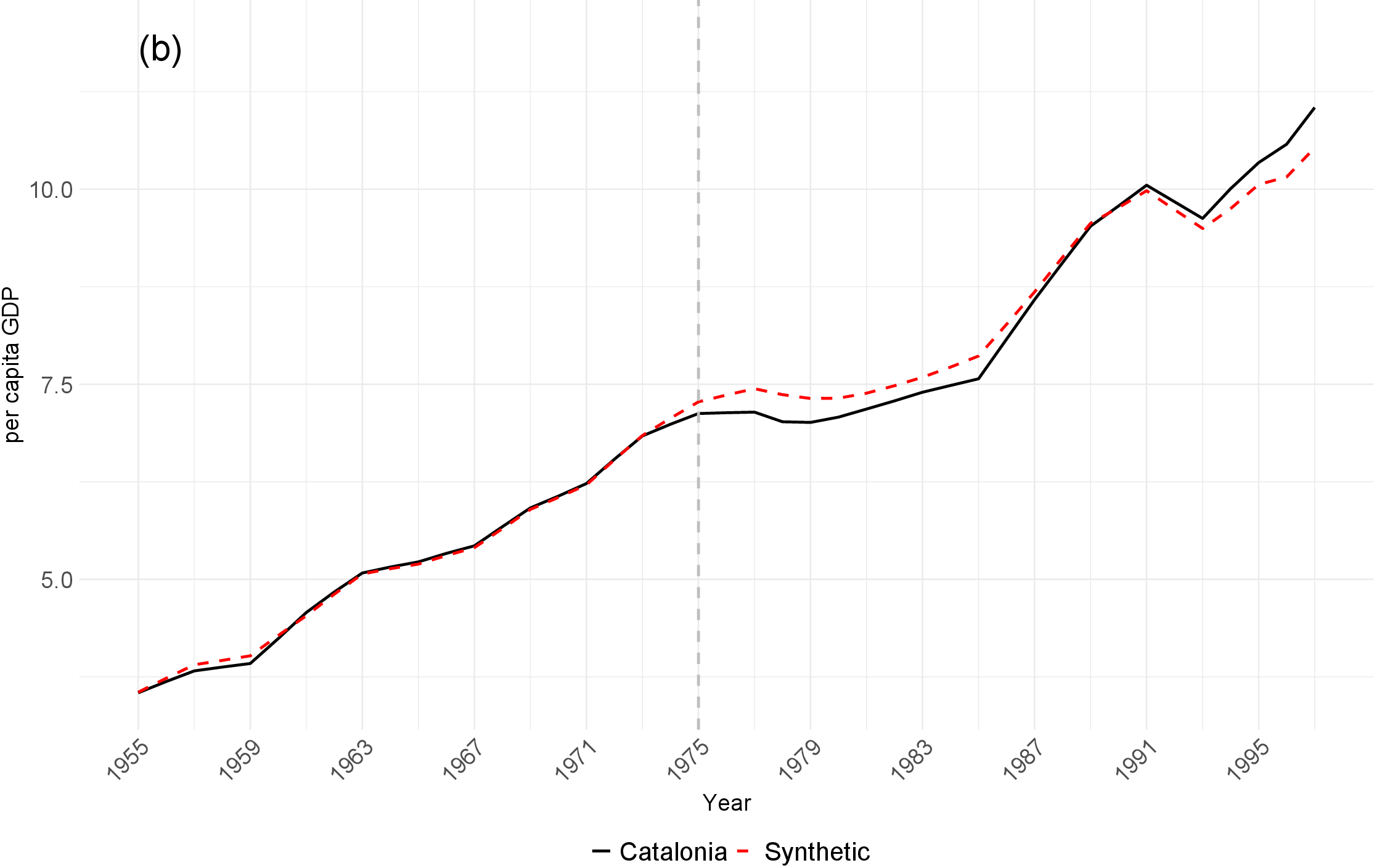}
\caption{Trends in per capita GDP: (a) actual Basque Country with terrorism versus the synthetic Basque Country obtained by {\sf RRMA}; (b) actual Catalonia without terrorism versus the synthetic Catalonia obtained by {\sf RRMA}. }
\label{Trends}
\label{figureA1}
\end{figure}

Figure \ref{figureA1} (a) shows that the synthetic \emph{Basque Country} obtained by {\sf RRMA} provides a statistical approximation to the actual per capita GDP of \emph{Basque Country} for 1955-1974 in the absence of terrorism. However, a natural question arises as to whether the observed gap between \emph{Basque Country} and its synthetic counterpart during the post-1975 period can be genuinely attributed to the impact of terrorism. To address this concern, we conduct a placebo study.
We choose \emph{Catalonia} (without terrorism) to carry out a placebo experiment because of its high similarity to \emph{Basque Country} in the pre-terrorism period regarding economic growth determinants.
Figure \ref{figureA1} (b) displays the per capita GDP for \emph{Catalonia} and
the ``synthetic \emph{Catalonia}" obtained by {\sf RRMA} during the period 1955--1997. Note that the ``synthetic \emph{Catalonia}" is constructed as a weighted combination of other 15 Spanish control regions. This figure clearly shows that the per capita GDP for \emph{Catalonia} can be reasonably reproduced by the proposed {\sf RRMA} over the entire period.
The placebo study concludes that the economic downturn in \emph{Basque Country} was mainly attributed to terrorist activities in the region and the proposed {\sf RRMA} can provide an accurate description for the unobserved per capita GDP of \emph{Basque Country} in 1975--1997 in the absence of terrorism.

\section*{Supplementary Material}
The supplementary material provides additional numerical results of the article and the technical proofs of lemmas, two Propositions and Theorems 1--6.

{\small
\setlength{\itemsep}{0.1em}
\setlength{\parsep}{0pt}
\bibliographystyle{apalike}
\bibliography{ref}
}

\end{document}